\documentclass[11pt, oneside]{article}   	
\usepackage{geometry}                		
\usepackage{graphicx}
\usepackage{float}				
\usepackage{amssymb}
\usepackage{amsmath}
\usepackage{authblk}
\usepackage{verbatim}

\usepackage{tikz}
\usetikzlibrary{shapes}
\tikzstyle{axis} = [line width = 0.8]
\tikzstyle{node} = [draw, rectangle, rounded corners, line width=0.5]
\tikzstyle{arrow} = [thick,->,>=stealth]

\usepackage[T1]{fontenc}
\usepackage[utf8]{inputenc}
\usepackage{tabularx,ragged2e,booktabs,caption}
\newcolumntype{C}[1]{>{\Centering}m{#1}}

\title{Fractals corresponding to the metallic means sequences}

\author{Y. S. J. Liang$^{1}$ \\  \texttt{yeishan@outlook.com} 
   \and  Darren C. Ong$^{1,2}$ \\  \texttt{darrenong@xmu.edu.my}}
\date{}							
\affil{\footnotesize{$^{1}$Department of Mathematics and Applied Mathematics, Xiamen University Malaysia, Jalan Sunsuria, Bandar Sunsuria, 43900, Sepang, Selangor, Malaysia \\ $^{2}$School of Mathematical Sciences, Xiamen University, Xiamen, 361005, Fujian, China}}

\begin{document}
\maketitle
\begin{abstract}
In this work, we introduce a class of fractal subsets of $[0,1]$ corresponding to the aperiodically ordered metallic means sequences. We find simple formulas for the fractal dimension for these fractals.

\end{abstract}

\section{Introduction}
In this paper, we consider fractal subsets of $[0,1]$ corresponding to a family of aperiodically ordered sequences known as the metallic means sequences. These are a family of infinite sequences of two symbols. A good reference for these sequences is \cite[Chapter 4]{Baake1999AGT}. From these sequences we obtain a tiling of $[0,1]$ with tiles of two different lengths. By selectively removing some tiles, and then re-tiling the remaining tiles the same way iteratively, we end up with a fractal set construction similar to a Cantor set. 

These are examples of dynamically generated Cantor sets considered in \cite{mori2001} and \cite{ichikawa2003}. We may also view them as $1$-dimensional analogues of the fractals considered in \cite{KahHeng2023}.

We introduce these fractals, and demonstrate simple formulae to calculate their fractal dimension. These formulae are not explicit, since they in general involve solutions of higher-degree polynomials. We note also that formulae for Hausdorff dimensions of the more general class of dynamically generated Cantor set already exists (for example in \cite{mori2001} and \cite{ichikawa2003}). Nevertheless, these simple formulae we demonstrate for the similarity and Hausdorff dimensions of these fractals shed a new perspective on these objects.

\section{Background}

\subsection{The similarity dimension of the $1/3$ Cantor set}
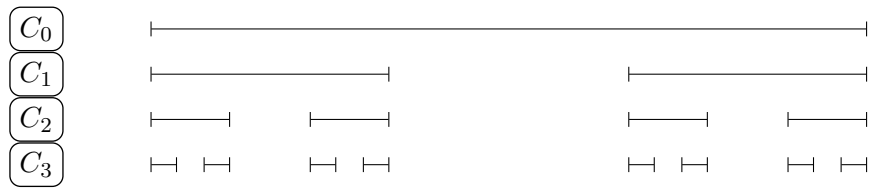
\begin{figure}[H]
\centering
\begin{tikzpicture}
	\node[node]at (-1.5,0.6) {$C_0$};	
	\draw[|-|] (0,0.6) -- (3*pi,0.6);
	
	\node[node]at (-1.5,0.0) {$C_1$};
	\draw[|-|] (0,0) -- (1*pi,0);    \draw[|-|] (2*pi,0) -- (3*pi,0);
	
	\node[node]at (-1.5,-0.6) {$C_2$};	
	\draw[|-|] (0,-.6) -- (1/3*pi,-0.6);\draw[|-|] (2/3*pi,-.6) -- (1*pi,-0.6);
	\draw[|-|] (2*pi,-.6) -- (7/3*pi,-0.6);\draw[|-|] (8/3*pi,-.6) -- (3*pi,-0.6);
	
	\node[node]at (-1.5,-1.2) {$C_3$};
	\draw[|-|] (0,-1.2) -- (1/9*pi,-1.2); \draw[|-|] (2/9*pi,-1.2) -- (3/9*pi,-1.2); 
	\draw[|-|] (6/9*pi,-1.2) -- (7/9*pi,-1.2); \draw[|-|] (8/9*pi,-1.2) -- (1*pi,-1.2);
	
	\draw[|-|] (2*pi,-1.2) -- (19/9*pi,-1.2); \draw[|-|] (20/9*pi,-1.2) -- (21/9*pi,-1.2);
	\draw[|-|] (24/9*pi,-1.2) -- (25/9*pi,-1.2); \draw[|-|] (26/9*pi,-1.2) -- (3*pi,-1.2);
	\end{tikzpicture}
\caption{Cantor set}
\label{canset}	
\end{figure}

Let us first consider the usual $1/3$ Cantor set (see for instance \cite{Strogatz94}). This $1/3$ Cantor set is the composed of two copies of itself, each scaled down by a factor of 3. This recursive construction is demonstrated below. The figures depicting $C_0$ to $C_3$ represent the first few iterations of the construction of this $1/3$ Cantor set in Figure \ref{canset}.

Let us review how to calculate the \emph{similarity dimension}, $d$ of the Cantor set. Recall that the similarity dimension is defined for a self-similar fractal, where the fractal is comprised of $m$ smaller copies of itself scaled down by a scale factor $r$. In this case, the similarity dimension $d$ is given by Equation \eqref{1}. 

\begin{equation}
\label{1}
\begin{aligned}[b]
d=\frac{\ln m}{\ln r} \\
\end{aligned}
\end{equation}

The $1/3$ Cantor set consists of two copies of itself scaled down by a scale factor of $3$, and so $m=2$ and $r=3$. The calculation for the similarity dimension is shown as below.

\begin{flalign*}
							d& =\log_{3}2  \\
							& \approx 0.6309
\end{flalign*}

\subsection{Hausdorff dimension of the $1/3$ Cantor set}

Let us quickly recall the calculation for the Hausdorff dimension of the Cantor set. For a subset $X$ of $[0,1]$, we define its $t$-dimensional Hausdorff measure by

$$H_t(X)=\inf_{\text {countable interval covers $\mathcal J$ of $X$}} \sum_{I\in\mathcal J} |I|^t
$$    

Then, the Hausdorff dimension is given by $dim_H(X)=\inf{\{t>0:H^t(X)=0\}}$.
 
Let $X$ be the usual $1/3$ Cantor Set from Figure \ref{canset}. We can cover the  Cantor set using intervals the usual way. For example, at the first stage, we have two intervals of length 1/3 and at the second stage, we have four intervals of length 1/9. Where $k$ corresponds to the $k$th stage of this Cantor set construction, $H^t(X)$ is the limit of the total length of intervals (raised to a power of $t$) when $k$ goes to infinity. Technically we have to consider all interval covers of $X$, not just the ones generated by the Cantor set contstruction. But we can argue that every interval in every possible cover of $X$ is contained within an interval from the  $k$th stage of a Cantor set construction, for some $k$ (for example, see Section 4.8 of \cite{Ziemer:2018te}. In other words
\begin{equation}
\label{tloi}
\begin{aligned}[b]
H_k^t(X) &=\lim_{k \to \infty} (number\ of\ intervals\ in\ step\ k)*(the\ length\ of\ the\ intervals\ in\ step\ k)^t
\end{aligned}
\end{equation}

Let us define 

\begin{equation}
\label{4b}
\begin{aligned}[b]
Y=(number\ of\ intervals\ in\ step\ 1)*(the\ length\ of\ the\ intervals\ in\ step\ 1)^t & \\
&
\end{aligned}
\end{equation}

and now notice that for the $k$th stage, $H^t_k(X)=Y^k$. And so

\begin{equation}
\label{4a}
\begin{aligned}[b]
H^t(X) &=\lim_{k \to \infty} {H_k^t(X)}  \in [0,\infty] \\
&=\lim_{k \to \infty} {Y^k}
\end{aligned}
\end{equation}

Thus to determine when $H^t(X)$ drops from infinity to 0, we have to calculate when $Y$ transitions from $Y>1$ to $Y<1$. Thus we can let $Y=1$ in Equation \ref{4b} to calculate the Hausdorff dimension in Equation \ref{thausdorffdc}. Let us label the number of intervals in stage 1 as $i$, and the length of the intervals in stage 1 as $j$. 

We then get from $Y=1$
\begin{equation}
\label{thausdorffdc}
\begin{aligned}[b]
ij^t=1\\
t=\frac{\ln \frac{1}{i}}{\ln j} \\
\end{aligned}
\end{equation}

The 1/3 Cantor set consists of 2 intervals with a length of $1/3$, and so $i=2$ and $j=1/3$. The calculation for the Hausdorff dimension is 

\begin{flalign*}
							& =\log_{\frac{1}{3}}\frac{1}{2}  \\
							& \approx 0.6309
\end{flalign*}


\subsection{Metallic means sequences and tilings of $[0,1]$}

The metallic means are defined as the positive root of the polynomial $x^2-px-q$ with $p,q \in N$. We know by Descartes' rule of signs (see for example \cite{wang2004}) that there is exactly one such positive root.

\captionof{table}{Some metallic means} \label{tab:title} 
    \begin{tabular}{ |p{3cm}|p{2cm}|p{2cm}| p{3cm}| p{3cm}| }
    \hline
     Metallic means & $p$ & $q$ & Symbol, $r$ & Value\\
    \hline
   Golden &1&1& $\phi$ & $\frac{1+\sqrt{5}}{2}$ \\
   Silver &2&1& $\delta$ & $1+\sqrt{2}$ \\
   Bronze &3&1& $\sigma$ & $\frac{3+\sqrt{13}}{2}$ \\
     Copper &1&2& $\alpha$ & $\frac{1+\sqrt{9}}{2} $\\
    Nickel &1&3& $\beta$ & $\frac{1+\sqrt{13}}{2}$ \\
  General &p&q& $\gamma$ & $\frac{\sqrt{p^2+4q}+p}{2}$ \\
    \hline
    \end{tabular}
    \label{table1}

The metallic means correspond to sequences of two letters in the following way. Consider an alphabet of two symbols $\mathcal A=\{a,b\}$. We then consider the free monoid $\mathcal A^*$ generated by this set. Then, we consider a mapping that takes $a$ to $a^pb^q$ and $b$ to $a$. There will be infinite sequences that serve as fixed points of this mapping. These are generated by starting with $a$ and applying the mapping to each letter. For example, in the case $p=q=1$ we have

$$a\to ab\to aba\to abaab\to abaababa \to abaababaabaab\to\ldots$$
The initial letters of this sequence of words stabilize, which generates the Fibonacci word, which is also the metallic mean sequence corresponding to the golden mean. See \cite[Section 4.1]{BaakeGrimm2013} for further details.The following table demonstrates the substitution rules for sequences corresponding to the other metallic means:

\


\section{The fractal corresponding to the Fibonacci Word}
\subsection{The Fibonacci tiling of $[0,1]$}
The first metallic means sequence is the the golden mean sequence, more commonly known as the \emph{Fibonacci word}. The sequence's  substitution rule, $\rho_F$ is

\begin{flalign*}
\rho_F : \ \ \ \ \ a &\rightarrow ab \\
b &\rightarrow a 
\end{flalign*}

Fibonacci word is a limit of a sequence of finite words of length $f_{i+2}$ where the $f_i $ is the $i$th Fibonacci number \cite{Leonardo1202}, where $f_0=0$ and $f_1=1$ \cite[Example 4.6]{BaakeGrimm2013}. \

Fibonacci word is generated by the addition of two previous items. As we can see from the substitution rule \cite[Example 4.3]{BaakeGrimm2013} on a finite alphabet $A=\{a,b\}$ with $a$ and $b$ is an endomorphism of the corresponding free group $F_n:=<a,b>$ generated by the letters of the alphabets \cite[Definition 4.1]{BaakeGrimm2013}. For every iterations, the $ab$ replaced the word $a$ and the $a$ replaced the word $b$. The examples for the first few finite words are $a\mapsto ab\mapsto  aba\mapsto  abaaba\mapsto  abaabaabaaba\mapsto ...$ and so on. 

The construction of the Fibonacci word corresponds to a tiling of the interval $[0,1]$ in the following way. We take the two letters $a,b$ of the alphabet to represent two types of tiles, where the $a$ represents a tile of length $\phi$ and $b$ represents a tile of length $1$. At each step of the Fibonacci word's construction, we scale down the length of both tiles by a factor of $\phi$. This series of tilings is depicted below in Figure  \ref{abb} \cite[Example 4.1]{BaakeGrimm2013}. 

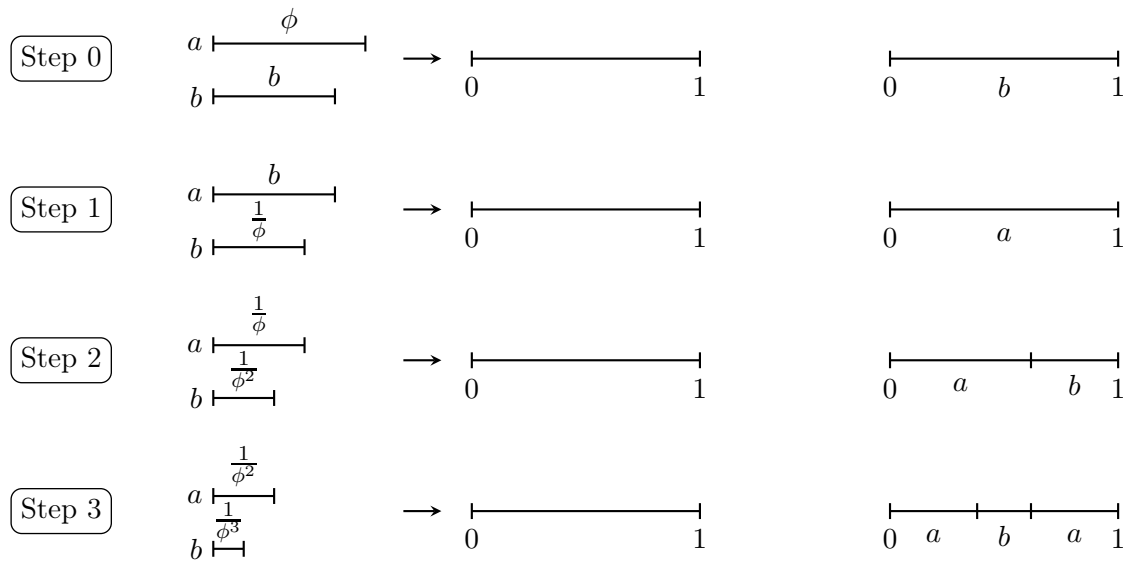
\begin{figure}[H]
\centering
\begin{tikzpicture}
		
		\draw[axis] (-4,1.7) -- (-2,1.7);
		\draw[axis] (-4,1.0) -- (-2.4,1.0);
		\draw[axis] (-0.6,1.5) -- (2.4,1.5);
		\draw[axis] (4.9,1.5) -- (7.9,1.5);
		\node[node]at (-6,1.5) {Step\ 0};
		
		\draw[shift={(-4,1.7)}, thick] (0,3pt) -- (0,-3pt);
		\path[shift={(-4,1.7)}, thick](1,3pt) -- (1,0pt) ;
		\path[shift={(-4,1.8)}, thick](0,3pt) -- (0,-3pt) node[left]{$a$};
		\path[shift={(-4,1.8)}, thick](1,3pt) -- (1,-3pt) node[above]{$\phi$};
		\draw[shift={(-4,1.7)}, thick] (2,3pt) -- (2,-3pt);
		
		\draw[shift={(-4,1.0)}, thick] (0,3pt) -- (0,-3pt);
		\path[shift={(-4,1.0)}, thick](1,3pt) -- (1,0pt) ;
		\path[shift={(-4,1.1)}, thick](0,3pt) -- (0,-3pt) node[left]{$b$};
		\path[shift={(-4,1.1)}, thick](0.8,3pt) -- (0.8,-3pt) node[above] {$b$};
		\draw[shift={(-4,1.0)}, thick] (1.6,3pt) -- (1.6,-3pt);
		
		\draw[arrow] (-1.5,1.5) --(-1,1.5);

		\draw[shift={(-4,1.5)}, thick] (3.4,3pt) -- (3.4,-3pt) node[below]{$0$};
		\draw[shift={(-4,1.5)}, thick] (6.4,3pt) -- (6.4,-3pt) node[below]{$1$};

		\draw[shift={(-4,1.5)}, thick] (8.9,3pt) -- (8.9,-3pt) node[below]{$0$};
		\path[shift={(-4,1.5)}, thick](10.4,3pt) -- (10.4,-3pt) node[below]{$b$};
		\draw[shift={(-4,1.5)}, thick] (11.9,3pt) -- (11.9,-3pt) node[below]{$1$};

		\draw[axis] (-4,-0.3) -- (-2.4,-0.3);
		\draw[axis] (-4,-1.0) -- (-2.8,-1.0);
		\draw[axis] (-0.6,-0.5) -- (2.4,-0.5);
		\draw[axis] (4.9,-0.5) -- (7.9,-0.5);
		\node[node]at (-6,-0.5) {Step\ 1};
		
		\draw[shift={(-4,-0.3)}, thick] (0,3pt) -- (0,-3pt);
		\path[shift={(-4,-0.3)}, thick](1,3pt) -- (1,0pt) ;
		\path[shift={(-4,-0.2)}, thick](0,3pt) -- (0,-3pt) node[left]{$a$};
		\path[shift={(-4,-0.2)}, thick](0.8,3pt) -- (0.8,-3pt) node[above] {$b$};
		\draw[shift={(-4,-0.3)}, thick] (1.6,3pt) -- (1.6,-3pt);
		
		\draw[shift={(-4,-1.0)}, thick] (0,3pt) -- (0,-3pt);
		\path[shift={(-4,-1.0)}, thick](1,3pt) -- (1,0pt) ;
		\path[shift={(-4,-0.9)}, thick](0,3pt) -- (0,-3pt) node[left]{$b$};
		\path[shift={(-4,-0.95)}, thick](0.6,3pt) -- (0.6,-3pt) node[above]{$\frac{1}{\phi}$};
		\draw[shift={(-4,-1.0)}, thick] (1.2,3pt) -- (1.2,-3pt);
		
		\draw[arrow] (-1.5,-0.5) --(-1,-0.5);

		\draw[shift={(-4,-0.5)}, thick] (3.4,3pt) -- (3.4,-3pt) node[below]{$0$};;
		\draw[shift={(-4,-0.5)}, thick] (6.4,3pt) -- (6.4,-3pt) node[below]{$1$};;

		\draw[shift={(-4,-0.5)}, thick] (8.9,3pt) -- (8.9,-3pt) node[below]{$0$};
		\path[shift={(-4,-0.5)}, thick](10.4,3pt) -- (10.4,-3pt) node[below]{$a$};
		\draw[shift={(-4,-0.5)}, thick] (11.9,3pt) -- (11.9,-3pt) node[below]{$1$};

		\draw[axis] (-4,-2.3) -- (-2.8,-2.3);
		\draw[axis] (-4,-3) -- (-3.2,-3);
		\draw[axis] (-0.6,-2.5) -- (2.4,-2.5);
		\draw[axis] (4.9,-2.5) -- (7.9,-2.5);
		\node[node]at (-6,-2.5) {Step\ 2};
		
		\draw[shift={(-4,-2.3)}, thick] (0,3pt) -- (0,-3pt);
		\path[shift={(-4,-2.3)}, thick](1,3pt) -- (1,0pt) ;
		\path[shift={(-4,-2.2)}, thick](0,3pt) -- (0,-3pt) node[left]{$a$};
		\path[shift={(-4,-2.2)}, thick](0.6,3pt) -- (0.6,-3pt) node[above]{$\frac{1}{\phi}$};
		\draw[shift={(-4,-2.3)}, thick] (1.2,3pt) -- (1.2,-3pt);
		
		\draw[shift={(-4,-3.0)}, thick] (0,3pt) -- (0,-3pt);
		\path[shift={(-4,-3.0)}, thick](1,3pt) -- (1,0pt) ;
		\path[shift={(-4,-2.9)}, thick](0,3pt) -- (0,-3pt) node[left]{$b$};
		\path[shift={(-4,-2.95)}, thick](0.4,3pt) -- (0.4,-3pt) node[above]{$\frac{1}{\phi^2}$};
		\draw[shift={(-4,-3)}, thick] (0.8,3pt) -- (0.8,-3pt);
		
		\draw[arrow] (-1.5,-2.5) --(-1,-2.5);

		\draw[shift={(-4,-2.5)}, thick] (3.4,3pt) -- (3.4,-3pt) node[below]{$0$};
		\draw[shift={(-4,-2.5)}, thick] (6.4,3pt) -- (6.4,-3pt) node[below]{$1$};

		\draw[shift={(-4,-2.5)}, thick] (8.9,3pt) -- (8.9,-3pt) node[below]{$0$};
		\path[shift={(-4,-2.5)}, thick](9.827,3pt) -- (9.827,-3pt) node[below]{$a$};
		\draw[shift={(-4,-2.5)}, thick] (10.754,3pt) -- (10.754,-3pt);
		\path[shift={(-4,-2.5)}, thick](10.754,3pt) -- (10.754,0pt) ;
		\path[shift={(-4,-2.45)}, thick](11.327,3pt) -- (11.327,-3pt) node[below]{$b$};
		\draw[shift={(-4,-2.5)}, thick] (11.9,3pt) -- (11.9,-3pt) node[below]{$1$};
		
		\draw[axis] (-4,-4.3) -- (-3.2,-4.3);
		\draw[axis] (-4,-5) -- (-3.6,-5);
		\draw[axis] (-0.6,-4.5) -- (2.4,-4.5);
		\draw[axis] (4.9,-4.5) -- (7.9,-4.5);
		\node[node]at (-6,-4.5) {Step\ 3};
		
		\draw[shift={(-4,-4.3)}, thick] (0,3pt) -- (0,-3pt);
		\path[shift={(-4,-4.3)}, thick](1,3pt) -- (1,0pt) ;
		\path[shift={(-4,-4.2)}, thick](0,3pt) -- (0,-3pt) node[left]{$a$};
		\path[shift={(-4,-4.2)}, thick](0.4,3pt) -- (0.4,-3pt) node[above]{$\frac{1}{\phi^2}$};
		\draw[shift={(-4,-4.3)}, thick] (0.8,3pt) -- (0.8,-3pt);
		
		\draw[shift={(-4,-5.0)}, thick] (0,3pt) -- (0,-3pt);
		\path[shift={(-4,-5.0)}, thick](1,3pt) -- (1,0pt) ;
		\path[shift={(-4,-4.9)}, thick](0,3pt) -- (0,-3pt) node[left]{$b$};
		\path[shift={(-4,-4.95)}, thick](0.2,3pt) -- (0.2,-3pt) node[above]{$\frac{1}{\phi^3}$};
		\draw[shift={(-4,-5)}, thick] (0.4,3pt) -- (0.4,-3pt);
		
		\draw[arrow] (-1.5,-4.5) --(-1,-4.5);

		\draw[shift={(-4,-4.5)}, thick] (3.4,3pt) -- (3.4,-3pt) node[below]{$0$};
		\draw[shift={(-4,-4.5)}, thick] (6.4,3pt) -- (6.4,-3pt) node[below]{$1$};

		\draw[shift={(-4,-4.5)}, thick] (8.9,3pt) -- (8.9,-3pt) node[below]{$0$};
		\path[shift={(-4,-4.5)}, thick](9.472955,3pt) -- (9.472955,-3pt) node[below]{$a$};
		\draw[shift={(-4,-4.5)}, thick] (10.04591,3pt) -- (10.04591,-3pt);
		\draw[shift={(-4,-4.5)}, thick](10.754,3pt) -- (10.754,-3pt) ;
		\path[shift={(-4,-4.45)}, thick](10.4,3pt) -- (10.4,-3pt) node[below]{$b$};
		\path[shift={(-4,-4.5)}, thick](11.327,3pt) -- (11.327,-3pt) node[below]{$a$};
		\draw[shift={(-4,-4.5)}, thick] (11.9,3pt) -- (11.9,-3pt) node[below]{$1$};

\end{tikzpicture}
\caption{Tilings corresponding to the Fibonacci inflation rule.}
\label{abb}
\end{figure}

Let us demonstrate that with this inflation factor, the sum of lengths of the tiles is $1$ at every step (so the tiling covers the interval $[0,1]$). We will show this by induction. For the base case, we simply observe that at the step $0$ case, we have a single tile of type $b$ which has length $1$. Assume now that the tiles of the step $k$ tiling have total length exactly $1$. We will show that the step $k+1$ tiling has the same total length. Consider a type $a$ tile in step $k$. It is replaced by a type $a$ tile and a type $b$ tile in step $k+1$, with these two tiles having a deflation factor $1/\phi$. So if $L$ is the length of the type $a$ tile in step $k$, the type $a$ tile in step $k+1$ has length $L/ \phi$, and the type $b$ tile in step $k+1$ has length $L/\phi^2$. The total length of these two tiles in step $k+1$ is $L/\phi+L/\phi^2=L$, which demonstrates that the $a$ tile in step $k$ gets replaced by two tiles of the same total length. On the other hand, a type $b$ tile in step $k$ gets replaced by a type $a$ tile in step $k+1$. By definition a type $a$ tile is $\phi$ times the length of the $b$ tile, so if we take into account the deflation factor of $\phi$, the type $a$ tile in step $k+1$ will have the same length as the type $b$ tile in Step $k$. Thus the transition from step $k$ to step $k+1$ does not change the total length of the tiles, which means at Step $k+1$ the total length of the tiles is still $1$. This concludes the induction proof.

We can derive in a constructive way that the appropriate scaling factor was $r=\phi$, the golden ratio. In order for the preceding proof to work, we call the scaling factor $r$ and the calculation is shown as below.

\begin{flalign*}
\frac{1}{r} + \frac{1}{r^2} 			&=1 \\ 
r^2-r-1						&=0 \\
r							&=\frac{1\pm \sqrt{5}}{2}. 	
\end{flalign*}

If we define a scaling factor to be larger than $1$, this means the inflation factor $r$ is equal to the golden ratio $\phi=\frac{1+\sqrt 5}{2}$.

\subsection{Constructing the fractal}

We will construct a family of fractal subsets of $[0,1]$, similar to the method in \cite{KahHeng2023}. The process proceeds as follows. For $n$ a positive integer, we construct a step $n$ tiling of $[0,1]$ as in Figure \ref{abb}. For example, when n=1 is a, when n=2 is ab, when n=3 is aba, when n=4 is abaab... .Then we remove $\alpha$ tiles of type $a$ and $\beta$ tiles of type $b$ from it. Then for each remaining tile in the tiling (both of type $a$ and $b$) we perform the same process. Each tile is itself tiled using a type $a$ tiling (with the tiles now scaled down appropriately) and we remove $\alpha$ tiles of type $a$ and $\beta$ tiles of type $b$. We will call the fractal generated this way the $(n,\alpha,\beta)$ fractal corresponding to the Fibonacci word. 

\begin{figure}[H]
\centering
\begin{tikzpicture}	
		\node[shift={(0,-3)}, node] at (-1,0) {tiling};
		\draw[shift={(0,-3)}, axis] (0,0) -- (2*pi,0);
		\draw[shift={(0,-3)}, thick] (0,3pt) -- (0,-3pt) node[below] {$0$};
		\path[shift={(0,-3)}] (0.382*pi,3pt) -- (0.382*pi,-3pt) node[below] {$\frac{1}{\phi^2}$};
		\path[shift={(0,-3)}] (0.382*pi,3pt) -- (0.382*pi,0pt) node[above] {$a$};
		\draw[shift={(0,-3)}, thick] (0.764*pi,3pt) -- (0.764*pi,-3pt);
		
		\path[shift={(0,-3)}] (1*pi,3pt) -- (1*pi,-3pt) node[below] {$\frac{1}{\phi^3}$};
		\path[shift={(0,-3)}] (1*pi,3pt) -- (1*pi,0pt) node[above] {$b$};
		\draw[shift={(0,-3)}, thick] (1.236*pi,3pt) -- (1.236*pi,-3pt);
		
		\path[shift={(0,-3)}] (1.618*pi,3pt) -- (1.618*pi,-3pt) node[below] {$\frac{1}{\phi^2}$};
		\path[shift={(0,-3)}] (1.618*pi,3pt) -- (1.618*pi,0pt) node[above] {$a$};
		\draw[shift={(0,-3)}, thick] (2*pi,3pt) -- (2*pi,-3pt) node[below] {$1$};
		
		\draw[arrow] (1*pi,-4.2) --(1*pi,-4.6);
		\node[shift={(0,-5.3)}, node] at (-1.5,0) {$1$st removal};
		\draw[shift={(0,-5.3)}, axis] (0,0) -- (0.764*pi,0);
		\draw[shift={(0,-5.3)}, axis] (1.236*pi,0) -- (2*pi,0);
		\draw[shift={(0,-5.3)}, thick] (0,3pt) -- (0,-3pt) node[below] {$0$};
		
		\path[shift={(0,-5.3)}] (0.382*pi,3pt) -- (0.382*pi,-3pt) node[below] {$\frac{1}{\phi^2}$};
		\path[shift={(0,-5.3)}] (0.382*pi,3pt) -- (0.382*pi,0pt) node[above] {$a$};
		\draw[shift={(0,-5.3)}, thick] (0.764*pi,3pt) -- (0.764*pi,-3pt);
		
		\draw[shift={(0,-5.3)}, thick] (1.236*pi,3pt) -- (1.236*pi,-3pt);
		
		\path[shift={(0,-5.3)}] (1.618*pi,3pt) -- (1.618*pi,-3pt) node[below] {$\frac{1}{\phi^2}$};
		\path[shift={(0,-5.3)}] (1.618*pi,3pt) -- (1.618*pi,0pt) node[above] {$a$};
		\draw[shift={(0,-5.3)}, thick] (2*pi,3pt) -- (2*pi,-3pt) node[below] {$1$};
		
		\draw[arrow] (1*pi,-6.5) --(1*pi,-6.9);
		\node[shift={(0,-7.6)}, node] at (-1.5,0) {$2$nd removal};
		\draw[shift={(0,-7.6)}, axis] (0,0) -- (0.2865*pi,0);
		\draw[shift={(0,-7.6)}, axis] (0.4775*pi,0) -- (0.764*pi,0);
		\draw[shift={(0,-7.6)}, axis] (1.236*pi,0) -- (1.5225*pi,0);
		\draw[shift={(0,-7.6)}, axis] (1.7135*pi,0) -- (2*pi,0);
		\draw[shift={(0,-7.6)}, thick] (0,3pt) -- (0,-3pt) node[below] {$0$};
		\path[shift={(0,-7.6)}] (0.1432*pi,3pt) -- (0.1432*pi,-3pt) node[below] {$\frac{1}{\phi^4}$};
		\path[shift={(0,-7.6)}] (0.1432*pi,3pt) -- (0.1432*pi,0pt) node[above] {$a$};
		
		\draw[shift={(0,-7.6)}, thick] (0.2865*pi,3pt) -- (0.2865*pi,-3pt);
		\draw[shift={(0,-7.6)}, thick] (0.4775*pi,3pt) -- (0.4775*pi,-3pt);
		
		\path[shift={(0,-7.6)}] (0.6207*pi,3pt) -- (0.6207*pi,-3pt) node[below] {$\frac{1}{\phi^4}$};
		\path[shift={(0,-7.6)}] (0.6207*pi,3pt) -- (0.6207*pi,0pt) node[above] {$a$};
		\draw[shift={(0,-7.6)}, thick] (0.764*pi,3pt) -- (0.764*pi,-3pt);
		
		\draw[shift={(0,-7.6)}, thick] (1.5225*pi,3pt) -- (1.5225*pi,-3pt);
		\draw[shift={(0,-7.6)}, thick] (1.7135*pi,3pt) -- (1.7135*pi,-3pt);
		
		\draw[shift={(0,-7.6)}, thick] (1.236*pi,3pt) -- (1.236*pi,-3pt);
		\path[shift={(0,-7.6)}] (1.3792*pi,3pt) -- (1.3792*pi,-3pt) node[below] {$\frac{1}{\phi^4}$};
		\path[shift={(0,-7.6)}] (1.3792*pi,3pt) -- (1.3792*pi,0pt) node[above] {$a$};
		
		\path[shift={(0,-7.6)}] (1.8567*pi,3pt) -- (1.8567*pi,-3pt) node[below] {$\frac{1}{\phi^4}$};
		\path[shift={(0,-7.6)}] (1.8567*pi,3pt) -- (1.8567*pi,0pt) node[above] {$a$};
		\draw[shift={(0,-7.6)}, thick] (2*pi,3pt) -- (2*pi,-3pt) node[below] {$1$};
		
		\draw[arrow] (1*pi,-8.8) --(1*pi,-9.2);
		\node[shift={(0,-9.9)}, node] at (-1.5,0) {$3$rd removal};
		\draw[shift={(0,-9.9)}, axis] (0,0) -- (0.1075*pi,0);
		\draw[shift={(0,-9.9)}, axis] (0.1791*pi,0) -- (0.2865*pi,0);
		\draw[shift={(0,-9.9)}, axis] (0.4775*pi,0) -- (0.585*pi,0);
		\draw[shift={(0,-9.9)}, axis] (0.6566*pi,0) -- (0.764*pi,0);
		\draw[shift={(0,-9.9)}, axis] (1.236*pi,0) -- (1.3435*pi,0);
		\draw[shift={(0,-9.9)}, axis] (1.4151*pi,0) -- (1.5225*pi,0);
		\draw[shift={(0,-9.9)}, axis] (1.7135*pi,0) -- (1.821*pi,0);
		\draw[shift={(0,-9.9)}, axis] (1.8926*pi,0) -- (2*pi,0);
		\draw[shift={(0,-9.9)}, thick] (0,3pt) -- (0,-3pt) node[below] {$0$};
		\draw[shift={(0,-9.9)}, thick] (0.1075*pi,3pt) -- (0.1075*pi,-3pt);
		
		\draw[shift={(0,-9.9)}, thick] (0.1791*pi,3pt) -- (0.1791*pi,-3pt);
		
		\draw[shift={(0,-9.9)}, thick] (0.2865*pi,3pt) -- (0.2865*pi,-3pt);
		\draw[shift={(0,-9.9)}, thick] (0.4775*pi,3pt) -- (0.4775*pi,-3pt);
		\draw[shift={(0,-9.9)}, thick] (0.585*pi,3pt) -- (0.585*pi,-3pt);
		
		\draw[shift={(0,-9.9)}, thick] (0.6566*pi,3pt) -- (0.6566*pi,-3pt);
		
		\draw[shift={(0,-9.9)}, thick] (1.3435*pi,3pt) -- (1.3435*pi,-3pt);
		
		\draw[shift={(0,-9.9)}, thick] (0.764*pi,3pt) -- (0.764*pi,-3pt);
		
		\draw[shift={(0,-9.9)}, thick] (1.236*pi,3pt) -- (1.236*pi,-3pt);
		\draw[shift={(0,-9.9)}, thick] (1.3435*pi,3pt) -- (1.3435*pi,-3pt);
		\draw[shift={(0,-9.9)}, thick] (1.4151*pi,3pt) -- (1.4151*pi,-3pt);
		
		\draw[shift={(0,-9.9)}, thick] (1.5225*pi,3pt) -- (1.5225*pi,-3pt);
		\draw[shift={(0,-9.9)}, thick] (1.7135*pi,3pt) -- (1.7135*pi,-3pt);
		\draw[shift={(0,-9.9)}, thick] (1.821*pi,3pt) -- (1.821*pi,-3pt);
		
		\draw[shift={(0,-9.9)}, thick] (1.8926*pi,3pt) -- (1.8926*pi,-3pt);
		
		\draw[shift={(0,-9.9)}, thick] (2*pi,3pt) -- (2*pi,-3pt) node[below] {$1$};
		
\end{tikzpicture}
\caption{Construction process for the $(3,0,1)$ fractal.}
\label{g3t}	
\end{figure}
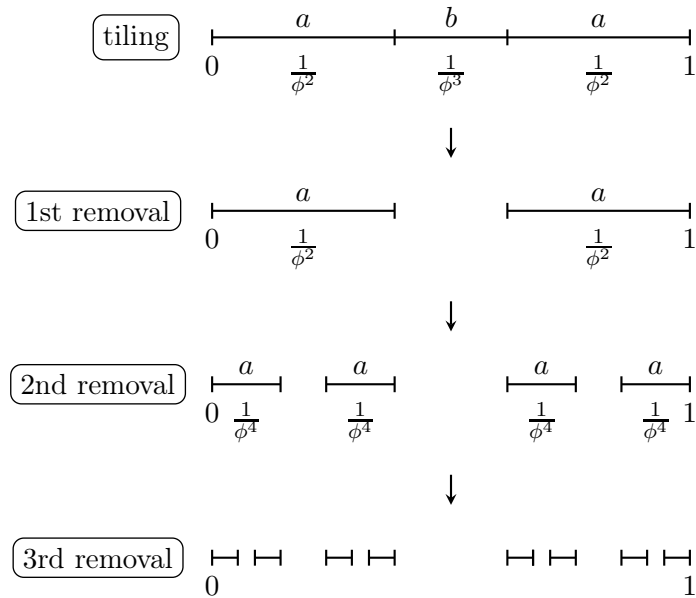

As an example,  Figure \ref{g3t} shows the fractal construction process for the $(3,0,1)$ fractal.

The calculation for the similarity dimension of this fractal is as follows

\begin{flalign*}
{(scaling\ factor, r)}^{(dimension, d)}   &= \text{total length of tiling after removals}   \\
(\phi^2)^d					    &=2	\\
d &=\frac{log 2}{log (\phi^2)} \\
d						&= 	0.7202
\end{flalign*}




Let us now calculate the Hausdorff dimension of this set, $dim_H(X)$. The method of calculating this Hausdorff dimension is not too different from how we calculated the Cantor Set. Recall that the definition of the Hausdorff dimension is $dim_H(X)=\inf{\{t>0:H^t(X)=0\}}$, where $H_t(X)$ refers to the t-dimensional Hausdorff measure of the set $X$.

As in our Cantor set calculation, it suffices to consider only interval covers that are generated from the steps of the fractal construction process. For example, in the first step we obtain a cover of $X$ consisting of two intervals of length $1/\phi^2$ each, from the second step we obtain a cover consisting of four intervals of length $1\phi^4$ each. In general, for the $k$th step our cover consists of $2^k$ intervals each of length $1/\phi^{2k}$. We label the $t$-dimensional Hausdorff measure of this $k$th-step cover as $H_k^t= 2^k/\phi^{2tk}$.  Then $H^t(X)=\lim_{k\to\infty} H_k^t(X)$.

Let us label $Y=2/\phi^{2t}$, then equivalently  $H^t(X)=\lim_{k\to\infty} Y^k$. 
Thus at $Y=1$, $H^t(X)$ drops from infinity to 0. We can then solve for the value of $t$ where this drop occurs, as follows.  \\

\begin{flalign*}
Y&=1 \\
2*\bigg(\frac{1}{\phi^2}\bigg)^t &=1 \\
\frac{2}{(\phi^{2t})} &=1 \\
2 &=\phi^{2t} \\
\ln 2 &=2t \ln \phi \\
\frac{\ln 2}{\ln \phi} &=2t \\
\frac{\ln 2}{2 \ln \phi} &=t \\
t &= 0.7202
\end{flalign*}

From the calculation, we know that $dim_H(X)=0.7202$.  \\

Let us look at a slightly more complicated example, the $(4,1,1)$ fractal. That is, we perform a $4$-step tiling as in Figure \ref{abb}, and then remove a tile of type $a$ and $b$ from it, and then recursively perform the same tiling and removal process to each remaining tile. This fractal construction process is visualized in Figure \ref{g42t} below:

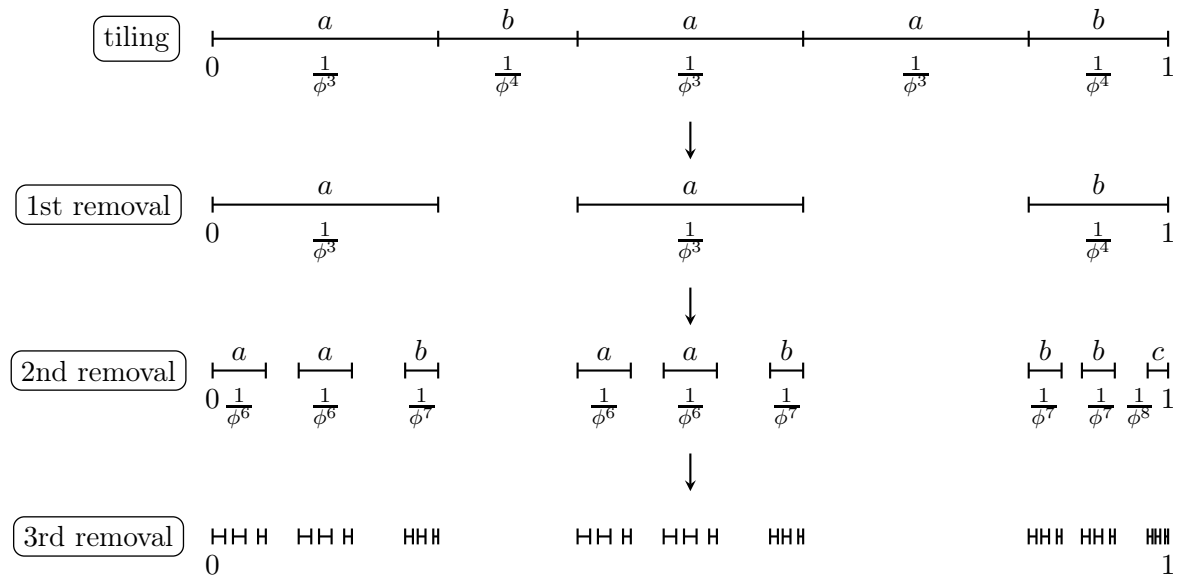
\begin{figure}[H]
\centering
\begin{tikzpicture}	
		\node[shift={(0,-4.5)}, node] at (-1,0) {tiling};
		\draw[shift={(0,-4.5)}, axis] (0,0) -- (4*pi,0);
		\draw[shift={(0,-4.5)}, thick] (0,3pt) -- (0,-3pt) node[below] {$0$};
		\path[shift={(0,-4.5)}] (0.11803*4*pi,3pt) -- (0.11803*4*pi,-3pt) node[below] {$\frac{1}{\phi^3}$};
		\path[shift={(0,-4.5)}] (0.11803*4*pi,3pt) -- (0.11803*4*pi,0pt) node[above] {$a$};
		\draw[shift={(0,-4.5)}, thick] (0.23607*4*pi,3pt) -- (0.23607*4*pi,-3pt);
		\path[shift={(0,-4.5)}] (0.30902*4*pi,3pt) -- (0.30902*4*pi,-3pt) node[below] {$\frac{1}{\phi^4}$};
		\path[shift={(0,-4.5)}] (0.30902*4*pi,3pt) -- (0.30902*4*pi,0pt) node[above] {$b$};
		\draw[shift={(0,-4.5)}, thick] (0.38197*4*pi,3pt) -- (0.38197*4*pi,-3pt);
		\path[shift={(0,-4.5)}] (0.5*4*pi,3pt) -- (0.5*4*pi,-3pt) node[below] {$\frac{1}{\phi^3}$};
		\path[shift={(0,-4.5)}] (0.5*4*pi,3pt) -- (0.5*4*pi,0pt) node[above] {$a$};
		\draw[shift={(0,-4.5)}, thick] (0.61804*4*pi,3pt) -- (0.61804*4*pi,-3pt);
		\path[shift={(0,-4.5)}] (0.73607*4*pi,3pt) -- (0.73607*4*pi,-3pt) node[below] {$\frac{1}{\phi^3}$};
		\path[shift={(0,-4.5)}] (0.73607*4*pi,3pt) -- (0.73607*4*pi,0pt) node[above] {$a$};
		\draw[shift={(0,-4.5)}, thick] (0.85411*4*pi,3pt) -- (0.85411*4*pi,-3pt);
		\path[shift={(0,-4.5)}] (0.92706*4*pi,3pt) -- (0.92706*4*pi,-3pt) node[below] {$\frac{1}{\phi^4}$};
		\path[shift={(0,-4.5)}] (0.92706*4*pi,3pt) -- (0.92706*4*pi,0pt) node[above] {$b$};
		\draw[shift={(0,-4.5)}, thick] (4*pi,3pt) -- (4*pi,-3pt) node[below] {$1$};

		\draw[arrow] (2*pi,-5.6) --(2*pi,-6.1);
		\node[shift={(0,-6.7)}, node] at (-1.5,0) {$1$st removal};
		\draw[shift={(0,-6.7)}, axis] (0,0) -- (0.23607*4*pi,0);
		\draw[shift={(0,-6.7)}, axis] (0.38197*4*pi,0) -- (0.61804*4*pi,0);
		\draw[shift={(0,-6.7)}, axis] (0.85411*4*pi,0) -- (4*pi,0);
		\draw[shift={(0,-6.7)}, thick] (0,3pt) -- (0,-3pt) node[below] {$0$};
		\path[shift={(0,-6.7)}] (0.11803*4*pi,3pt) -- (0.11803*4*pi,-3pt) node[below] {$\frac{1}{\phi^3}$};
		\path[shift={(0,-6.7)}] (0.11803*4*pi,3pt) -- (0.11803*4*pi,0pt) node[above] {$a$};
		\draw[shift={(0,-6.7)}, thick] (0.23607*4*pi,3pt) -- (0.23607*4*pi,-3pt);
		
		\draw[shift={(0,-6.7)}, thick] (0.38197*4*pi,3pt) -- (0.38197*4*pi,-3pt);
		\path[shift={(0,-6.7)}] (0.5*4*pi,3pt) -- (0.5*4*pi,-3pt) node[below] {$\frac{1}{\phi^3}$};
		\path[shift={(0,-6.7)}] (0.5*4*pi,3pt) -- (0.5*4*pi,0pt) node[above] {$a$};
		\draw[shift={(0,-6.7)}, thick] (0.61804*4*pi,3pt) -- (0.61804*4*pi,-3pt);
		
		\draw[shift={(0,-6.7)}, thick] (0.85411*4*pi,3pt) -- (0.85411*4*pi,-3pt);
		\path[shift={(0,-6.7)}] (0.92706*4*pi,3pt) -- (0.92706*4*pi,-3pt) node[below] {$\frac{1}{\phi^4}$};
		\path[shift={(0,-6.7)}] (0.92706*4*pi,3pt) -- (0.92706*4*pi,0pt) node[above] {$b$};
		\draw[shift={(0,-6.7)}, thick] (4*pi,3pt) -- (4*pi,-3pt) node[below] {$1$};

		\draw[arrow] (2*pi,-7.8) --(2*pi,-8.3);
		\node[shift={(0,-8.9)}, node] at (-1.5,0) {$2$nd removal};
		\draw[shift={(0,-8.9)}, axis] (0,0) -- (0.0557*4*pi,0);
		\draw[shift={(0,-8.9)}, axis] (0.0901*4*pi,0) -- (0.1458*4*pi,0);
		\draw[shift={(0,-8.9)}, axis] (0.2015*4*pi,0) -- (0.236*4*pi,0);
		
		\draw[shift={(0,-8.9)}, axis] (0.38197*4*pi,0) -- (0.43767*4*pi,0);
		\draw[shift={(0,-8.9)}, axis] (0.47207*4*pi,0) -- (0.52777*4*pi,0);
		\draw[shift={(0,-8.9)}, axis] (0.58347*4*pi,0) -- (0.618*4*pi,0);

		\draw[shift={(0,-8.9)}, axis] (0.85411*4*pi,0) -- (0.88851*4*pi,0);
		\draw[shift={(0,-8.9)}, axis] (0.90981*4*pi,0) -- (0.94421*4*pi,0);
		\draw[shift={(0,-8.9)}, axis] (0.97861*4*pi,0) -- (4*pi,0);

		\draw[shift={(0,-8.9)}, thick] (0,3pt) -- (0,-3pt) node[below] {$0$};
		
		\draw[shift={(0,-8.9)}, thick] (0.0557*4*pi,3pt) -- (0.0557*4*pi,-3pt);
		\path[shift={(0,-8.9)}] (0.02785*4*pi,3pt) -- (0.02785*4*pi,-3pt) node[below] {$\frac{1}{\phi^6}$};
		\path[shift={(0,-8.9)}] (0.02785*4*pi,3pt) -- (0.02785*4*pi,0pt) node[above] {$a$};
		\draw[shift={(0,-8.9)}, thick] (0.0901*4*pi,3pt) -- (0.0901*4*pi,-3pt);
		\draw[shift={(0,-8.9)}, thick] (0.1458*4*pi,3pt) -- (0.1458*4*pi,-3pt);
		\draw[shift={(0,-8.9)}, thick] (0.2015*4*pi,3pt) -- (0.2015*4*pi,-3pt);
		\path[shift={(0,-8.9)}] (0.21785*4*pi,3pt) -- (0.21875*4*pi,-3pt) node[below] {$\frac{1}{\phi^7}$};
		\path[shift={(0,-8.9)}] (0.21785*4*pi,3pt) -- (0.21785*4*pi,0pt) node[above] {$b$};
		
		\path[shift={(0,-8.9)}] (0.11803*4*pi,3pt) -- (0.11803*4*pi,-3pt) node[below] {$\frac{1}{\phi^6}$};
		\path[shift={(0,-8.9)}] (0.11803*4*pi,3pt) -- (0.11803*4*pi,0pt) node[above] {$a$};
		\draw[shift={(0,-8.9)}, thick] (0.23607*4*pi,3pt) -- (0.23607*4*pi,-3pt);
		
		\draw[shift={(0,-8.9)}, thick] (0.38197*4*pi,3pt) -- (0.38197*4*pi,-3pt);
		\path[shift={(0,-8.9)}] (0.40982*4*pi,3pt) -- (0.40982*4*pi,-3pt) node[below] {$\frac{1}{\phi^6}$};
		\path[shift={(0,-8.9)}] (0.40982*4*pi,3pt) -- (0.40982*4*pi,0pt) node[above] {$a$};
		
		\draw[shift={(0,-8.9)}, thick] (0.43767*4*pi,3pt) -- (0.43767*4*pi,-3pt);
		\draw[shift={(0,-8.9)}, thick] (0.47207*4*pi,3pt) -- (0.47207*4*pi,-3pt);
		\draw[shift={(0,-8.9)}, thick] (0.52777*4*pi,3pt) -- (0.52777*4*pi,-3pt);
		\draw[shift={(0,-8.9)}, thick] (0.58347*4*pi,3pt) -- (0.58347*4*pi,-3pt);
		\path[shift={(0,-8.9)}] (0.60067*4*pi,3pt) -- (0.60067*4*pi,-3pt) node[below] {$\frac{1}{\phi^7}$};
		\path[shift={(0,-8.9)}] (0.60067*4*pi,3pt) -- (0.60067*4*pi,0pt) node[above] {$b$};
		
		\path[shift={(0,-8.9)}] (0.5*4*pi,3pt) -- (0.5*4*pi,-3pt) node[below] {$\frac{1}{\phi^6}$};
		\path[shift={(0,-8.9)}] (0.5*4*pi,3pt) -- (0.5*4*pi,0pt) node[above] {$a$};
		\draw[shift={(0,-8.9)}, thick] (0.61804*4*pi,3pt) -- (0.61804*4*pi,-3pt);
		
		\draw[shift={(0,-8.9)}, thick] (0.85411*4*pi,3pt) -- (0.85411*4*pi,-3pt);
		\path[shift={(0,-8.9)}] (0.87*4*pi,3pt) -- (0.87*4*pi,-3pt) node[below] {$\frac{1}{\phi^{7}}$};
		\path[shift={(0,-8.9)}] (0.87131*4*pi,3pt) -- (0.87131*4*pi,0pt) node[above] {$b$};
		
		\draw[shift={(0,-8.9)}, thick] (0.88851*4*pi,3pt) -- (0.88851*4*pi,-3pt);
		\draw[shift={(0,-8.9)}, thick] (0.90981*4*pi,3pt) -- (0.90981*4*pi,-3pt);
		\draw[shift={(0,-8.9)}, thick] (0.94421*4*pi,3pt) -- (0.94421*4*pi,-3pt);
		\draw[shift={(0,-8.9)}, thick] (0.97861*4*pi,3pt) -- (0.97861*4*pi,-3pt);
		\path[shift={(0,-8.9)}] (0.97*4*pi,3pt) -- (0.97*4*pi,-3pt) node[below] {$\frac{1}{\phi^{8}}$};
		\path[shift={(0,-8.9)}] (0.98926*4*pi,3pt) -- (0.98926*4*pi,0pt) node[above] {$c$};
		
		\path[shift={(0,-8.9)}] (0.93*4*pi,3pt) -- (0.93*4*pi,-3pt) node[below] {$\frac{1}{\phi^{7}}$};
		\path[shift={(0,-8.9)}] (0.92706*4*pi,3pt) -- (0.92706*4*pi,0pt) node[above] {$b$};
		
		\draw[shift={(0,-8.9)}, thick] (4*pi,3pt) -- (4*pi,-3pt) node[below] {$1$};

		\draw[arrow] (2*pi,-10) --(2*pi,-10.5);
		\node[shift={(0,-11.1)}, node] at (-1.5,0) {$3$rd removal};
		\draw[shift={(0,-11.1)}, axis] (0,0) -- (0.01315*4*pi,0);
		\draw[shift={(0,-11.1)}, axis] (0.02128*4*pi,0) -- (0.03443*4*pi,0);

		\draw[shift={(0,-11.1)}, axis] (0.04758*4*pi,0) -- (0.05571*4*pi,0);
		
		\draw[shift={(0,-11.1)}, axis] (0.0901*4*pi,0) -- (0.10325*4*pi,0);
		\draw[shift={(0,-11.1)}, axis] (0.11138*4*pi,0) -- (0.12453*4*pi,0);
		\draw[shift={(0,-11.1)}, axis] (0.13768*4*pi,0) -- (0.14581*4*pi,0);

		\draw[shift={(0,-11.1)}, axis] (0.2015*4*pi,0) -- (0.20963*4*pi,0);
		\draw[shift={(0,-11.1)}, axis] (0.21466*4*pi,0) -- (0.22279*4*pi,0);
		\draw[shift={(0,-11.1)}, axis] (0.23092*4*pi,0) -- (0.236*4*pi,0);

		\draw[shift={(0,-11.1)}, axis] (0.38197*4*pi,0) -- (0.39512*4*pi,0);
		\draw[shift={(0,-11.1)}, axis] (0.40325*4*pi,0) -- (0.4164*4*pi,0);
		\draw[shift={(0,-11.1)}, axis] (0.42955*4*pi,0) -- (0.43767*4*pi,0);
		
		\draw[shift={(0,-11.1)}, axis] (0.47207*4*pi,0) -- (0.48522*4*pi,0);
		\draw[shift={(0,-11.1)}, axis] (0.49335*4*pi,0) -- (0.5065*4*pi,0);
		\draw[shift={(0,-11.1)}, axis] (0.51965*4*pi,0) -- (0.52777*4*pi,0);

		\draw[shift={(0,-11.1)}, axis] (0.58347*4*pi,0) -- (0.5916*4*pi,0);
		\draw[shift={(0,-11.1)}, axis] (0.59663*4*pi,0) -- (0.60476*4*pi,0);
		\draw[shift={(0,-11.1)}, axis] (0.61289*4*pi,0) -- (0.618*4*pi,0);

		\draw[shift={(0,-11.1)}, axis] (0.85411*4*pi,0) -- (0.86224*4*pi,0);
		\draw[shift={(0,-11.1)}, axis] (0.86726*4*pi,0) -- (0.87539*4*pi,0);
		\draw[shift={(0,-11.1)}, axis] (0.88352*4*pi,0) -- (0.88851*4*pi,0);
		
		\draw[shift={(0,-11.1)}, axis] (0.90981*4*pi,0) -- (0.91794*4*pi,0);
		\draw[shift={(0,-11.1)}, axis] (0.92296*4*pi,0) -- (0.93109*4*pi,0);
		\draw[shift={(0,-11.1)}, axis] (0.93922*4*pi,0) -- (0.94421*4*pi,0);

		\draw[shift={(0,-11.1)}, axis] (0.97861*4*pi,0) -- (0.98364*4*pi,0);
		\draw[shift={(0,-11.1)}, axis] (0.98676*4*pi,0) -- (0.99179*4*pi,0);
		\draw[shift={(0,-11.1)}, axis] (0.99682*4*pi,0) -- (4*pi,0);
		
		\draw[shift={(0,-11.1)}, thick] (0,3pt) -- (0,-3pt) node[below] {$0$};
		
		\draw[shift={(0,-11.1)}, thick] (0.01315*4*pi,3pt) -- (0.01315*4*pi,-3pt);
		\draw[shift={(0,-11.1)}, thick] (0.02128*4*pi,3pt) -- (0.02128*4*pi,-3pt);
		\draw[shift={(0,-11.1)}, thick] (0.03443*4*pi,3pt) -- (0.03443*4*pi,-3pt);
		\draw[shift={(0,-11.1)}, thick] (0.04758*4*pi,3pt) -- (0.04758*4*pi,-3pt);
		
		\draw[shift={(0,-11.1)}, thick] (0.0557*4*pi,3pt) -- (0.0557*4*pi,-3pt);
		\draw[shift={(0,-11.1)}, thick] (0.0901*4*pi,3pt) -- (0.0901*4*pi,-3pt);
		
		\draw[shift={(0,-11.1)}, thick] (0.10325*4*pi,3pt) -- (0.10325*4*pi,-3pt);
		\draw[shift={(0,-11.1)}, thick] (0.11138*4*pi,3pt) -- (0.11138*4*pi,-3pt);
		\draw[shift={(0,-11.1)}, thick] (0.12453*4*pi,3pt) -- (0.12453*4*pi,-3pt);
		\draw[shift={(0,-11.1)}, thick] (0.13768*4*pi,3pt) -- (0.13768*4*pi,-3pt);
		
		\draw[shift={(0,-11.1)}, thick] (0.1458*4*pi,3pt) -- (0.1458*4*pi,-3pt);
		\draw[shift={(0,-11.1)}, thick] (0.2015*4*pi,3pt) -- (0.2015*4*pi,-3pt);
		
		\draw[shift={(0,-11.1)}, thick] (0.20963*4*pi,3pt) -- (0.20963*4*pi,-3pt);
		\draw[shift={(0,-11.1)}, thick] (0.21466*4*pi,3pt) -- (0.21466*4*pi,-3pt);
		\draw[shift={(0,-11.1)}, thick] (0.22279*4*pi,3pt) -- (0.22279*4*pi,-3pt);
		\draw[shift={(0,-11.1)}, thick] (0.23092*4*pi,3pt) -- (0.23092*4*pi,-3pt);
		
		\draw[shift={(0,-11.1)}, thick] (0.23607*4*pi,3pt) -- (0.23607*4*pi,-3pt);
		
		\draw[shift={(0,-11.1)}, thick] (0.38197*4*pi,3pt) -- (0.38197*4*pi,-3pt);
		
		
		
		\draw[shift={(0,-11.1)}, thick] (0.39512*4*pi,3pt) -- (0.39512*4*pi,-3pt);
		\draw[shift={(0,-11.1)}, thick] (0.40325*4*pi,3pt) -- (0.40325*4*pi,-3pt);
		\draw[shift={(0,-11.1)}, thick] (0.4164*4*pi,3pt) -- (0.4164*4*pi,-3pt);
		\draw[shift={(0,-11.1)}, thick] (0.42955*4*pi,3pt) -- (0.42955*4*pi,-3pt);
		
		\draw[shift={(0,-11.1)}, thick] (0.43767*4*pi,3pt) -- (0.43767*4*pi,-3pt);
		\draw[shift={(0,-11.1)}, thick] (0.47207*4*pi,3pt) -- (0.47207*4*pi,-3pt);
		
		\draw[shift={(0,-11.1)}, thick] (0.48522*4*pi,3pt) -- (0.48522*4*pi,-3pt);
		\draw[shift={(0,-11.1)}, thick] (0.49335*4*pi,3pt) -- (0.49335*4*pi,-3pt);
		\draw[shift={(0,-11.1)}, thick] (0.5065*4*pi,3pt) -- (0.5065*4*pi,-3pt);
		\draw[shift={(0,-11.1)}, thick] (0.51965*4*pi,3pt) -- (0.51965*4*pi,-3pt);
		
		\draw[shift={(0,-11.1)}, thick] (0.52777*4*pi,3pt) -- (0.52777*4*pi,-3pt);
		\draw[shift={(0,-11.1)}, thick] (0.58347*4*pi,3pt) -- (0.58347*4*pi,-3pt);
		
		\draw[shift={(0,-11.1)}, thick] (0.5916*4*pi,3pt) -- (0.5916*4*pi,-3pt);
		\draw[shift={(0,-11.1)}, thick] (0.59663*4*pi,3pt) -- (0.59663*4*pi,-3pt);
		\draw[shift={(0,-11.1)}, thick] (0.60476*4*pi,3pt) -- (0.60476*4*pi,-3pt);
		\draw[shift={(0,-11.1)}, thick] (0.61289*4*pi,3pt) -- (0.61289*4*pi,-3pt);
		
		\draw[shift={(0,-11.1)}, thick] (0.61804*4*pi,3pt) -- (0.61804*4*pi,-3pt);
		
		\draw[shift={(0,-11.1)}, thick] (0.85411*4*pi,3pt) -- (0.85411*4*pi,-3pt);
		
		\draw[shift={(0,-11.1)}, thick] (0.86224*4*pi,3pt) -- (0.86224*4*pi,-3pt);
		\draw[shift={(0,-11.1)}, thick] (0.86726*4*pi,3pt) -- (0.86726*4*pi,-3pt);
		\draw[shift={(0,-11.1)}, thick] (0.87539*4*pi,3pt) -- (0.87539*4*pi,-3pt);
		\draw[shift={(0,-11.1)}, thick] (0.88352*4*pi,3pt) -- (0.88352*4*pi,-3pt);
		
		\draw[shift={(0,-11.1)}, thick] (0.88851*4*pi,3pt) -- (0.88851*4*pi,-3pt);
		\draw[shift={(0,-11.1)}, thick] (0.90981*4*pi,3pt) -- (0.90981*4*pi,-3pt);
		
		\draw[shift={(0,-11.1)}, thick] (0.91794*4*pi,3pt) -- (0.91794*4*pi,-3pt);
		\draw[shift={(0,-11.1)}, thick] (0.92296*4*pi,3pt) -- (0.92296*4*pi,-3pt);
		\draw[shift={(0,-11.1)}, thick] (0.93109*4*pi,3pt) -- (0.93109*4*pi,-3pt);
		\draw[shift={(0,-11.1)}, thick] (0.93922*4*pi,3pt) -- (0.93922*4*pi,-3pt);
		
		\draw[shift={(0,-11.1)}, thick] (0.94421*4*pi,3pt) -- (0.94421*4*pi,-3pt);
		\draw[shift={(0,-11.1)}, thick] (0.97861*4*pi,3pt) -- (0.97861*4*pi,-3pt);
		
		\draw[shift={(0,-11.1)}, thick] (0.98364*4*pi,3pt) -- (0.98364*4*pi,-3pt);
		\draw[shift={(0,-11.1)}, thick] (0.98676*4*pi,3pt) -- (0.98676*4*pi,-3pt);
		\draw[shift={(0,-11.1)}, thick] (0.99179*4*pi,3pt) -- (0.99179*4*pi,-3pt);
		\draw[shift={(0,-11.1)}, thick] (0.99682*4*pi,3pt) -- (0.99682*4*pi,-3pt);
		
		\draw[shift={(0,-11.1)}, thick] (4*pi,3pt) -- (4*pi,-3pt) node[below] {$1$};
		
\end{tikzpicture}
\caption{Construction of the $(4,1,1)$ fractal}
\label{g42t}	
\end{figure}

We calculate the similarity dimension, $d$ for the $(4,1,1)$ in a similar way as before: 

\begin{flalign*}
{(scaling\ factor, r)}^{(dimension, d)}   &=length\ of\ 4-tile'' \\
(\phi^d)						&=2(\frac{1}{\phi^3})+(\frac{1}{\phi^3}*\frac{1}{\phi^d}) \\
							&=\frac{2}{\phi^3}+\frac{1}{\phi^{3d}} \\
(\phi^4)^d						&=2 \phi^d+1 \\
Let\ \phi^d &=x 					\\
						x^4	&=2x+1 \\
						0	&=x^4-2x-1 \\
						x	& \approx 1.3953 \\
so\ \phi^d						&\approx 1.3953 \\
d							&\approx \log_{\phi}1.3953 \\
							& \approx 	0.6922	
\end{flalign*}

Let us now calculate the Hausdorff dimension. As before, we only consider interval covers that arise from the construction of the fractal. We can see that the sum of the lengths of intervals in the cover corresponding to the $k$th step of the construction is given by $(\frac{2}{\phi^3}+\frac{1}{\phi^4})^k$.
For example, when $k=1$ we obtain $
(\frac{2}{\phi^3}+\frac{1}{\phi^4}).$ We read this as $2$ intervals length $(\frac{1}{\phi^3})$ and $1$ interval length $(\frac{1}{\phi^4})$. 

When $k=2, 
(\frac{2}{\phi^3}+\frac{1}{\phi^4})^2 =\frac{4}{\phi^6}+\frac{4}{\phi^7}+\frac{1}{\phi^8}$, which corresponds to $4$ intervals length $(\frac{1}{\phi^6}), 4$ intervals length $(\frac{1}{\phi^7}), 1$ interval length $(\frac{1}{\phi^8})$.  In general, 
$(\frac{2}{\phi^3}+\frac{1}{\phi^4})^k=\frac{c_0}{\phi^{3k}}+\frac{c_1}{\phi^{3k+1}}+\frac{c_2}{\phi^{3k+2}}+...+\frac{c_k}{\phi^{4k}}$ which refers to $c_0$ intervals length $(\frac{1}{\phi^{3k}}), c_1$ intervals length $(\frac{1}{\phi^{3k+1}})... , c_k$ intervals length $(\frac{1}{\phi^{4k}})$. \\

In the first removal of  the(4, 1, 1) fractal as shown in Figure \ref{g42t}, there were 2 tiles of length $\frac{1}{\phi^3}$ labeled as $a$ and 1 tile of length $\frac{1}{\phi^4}$ labeled as $b$ because we removed 1 tile of length $\frac{1}{\phi^3}$ labeled as $a$ and 1 tile of length $\frac{1}{\phi^4}$ labeled as $b$. Then $H_k^t(X)$ as before is the sum of $\ (number\ of\ tiles\ at\ step\ k)*(length\ of\ tiles\ at\ step\ k)^t$ which is $(\frac{2}{\phi^3}+\frac{1}{\phi^4})^k$. And since $H_k^t(X)$ also same as $Y^k$, this means $Y=\frac{2}{\phi^{3t}}+\frac{1}{\phi^{4t}}$. In order to find the Hausdorff dimension, we let $Y$ be $1$ when $H^t(X)$ goes from infinity to $0$. \\

\begin{flalign*}
Y   &=1\\
(\frac{2}{\phi^{3t}}+\frac{1}{\phi^{4t}}) 	&=1 \\
2\phi^t+1						&=\phi^{4t} \\
Let\ \phi^t &=x \\
2x-x^4						&=-1 \\
x   							&\approx 1.3953 \\
\phi^t						&\approx 1.3953 \\
t							&\approx \log_{\phi}1.3953 \\
							&\approx 	0.6922	
\end{flalign*}

From the calculation, we know that $t$ for $(4,1,1)$ fractal is 0.6922. \\

For a general $(n,l,s)$ fractal generated from the tiling corresponding to the Fibonacci word, the dimension might not be expressible using elementary functions, since it is a root of a high degree polynomial. Explicitly,

\begin{flalign*}
{(scaling\ factor, r)}^{(dimension, d)}   =&length\ after\ substitution,\ l_a\  \\
(\phi^d)					=&(total\ length\ of\ long\ tile,\ t_{l})+(total\ length\ of\ short\ tile,\ t_{s} )*\phi^{-d} \\
							=&(a_{n}-l)+(a_{n-1}-s)*\phi^{-d} \\
(\phi^d)^{n}					=&(a_{n}-l)*\phi^d+(a_{n-1}-s) 
\end{flalign*}

Let $x=\tilde{x}$ be the positive solution (by \cite{wang2004} there is exactly one) to 
\begin{flalign*} 
 x^n  =&\bigg(\frac{1}{\sqrt{5}}\bigg(\frac{1}{2}+\frac{\sqrt{5}}{2}\bigg)^{n}-\frac{1}{\sqrt{5}}\bigg(\frac{1}{2}-\frac{\sqrt{5}}{2}\bigg)^{n}-l \bigg)*x+ \\ 
& \bigg(\frac{1}{\sqrt{5}}\bigg(\frac{1}{2}+\frac{\sqrt{5}}{2}\bigg)^{n-1}-\frac{1}{\sqrt{5}}\bigg(\frac{1}{2}-\frac{\sqrt{5}}{2}\bigg)^{n-1}-s \bigg) 	\\	 
\end{flalign*}
Then $\phi^d=\tilde  x$ and so	$d=\log_\phi \tilde x$. \\

For Hausdorff dimension in a general $(n,l,s)$ fractal:
\begin{flalign*}
H^t(x)=&(number\ of\ intervals)(the\ length\ of\ the\ intervals)^t \\
1=&(total\ length\ of\ long\ tiles\ intervals, t_l)+(total\ length\ of\ the\ short\ tile\ intervals, t_s) \\
1=&(a_{n}-l)(1/\phi^3)^t+(a_{n-1}-s)(1/\phi^4)^t \\
1=&(1/\phi^t)^4[(a_{n}-l)\phi^t+(a_{n-1}-s)] \\
\phi^{4t}=&(a_{n}-l)\phi^t+(a_{n-1}-s) \\
(\phi^t)^n=&(a_{n}-l)\phi^t+(a_{n-1}-s) 
\end{flalign*}

Let $x=\tilde{x}$ be the positive solution (again by \cite{wang2004} there is only one) to 
\begin{flalign*} 
 x^n  =&\bigg(\frac{1}{\sqrt{5}}\bigg(\frac{1}{2}+\frac{\sqrt{5}}{2}\bigg)^{n}-\frac{1}{\sqrt{5}}\bigg(\frac{1}{2}-\frac{\sqrt{5}}{2}\bigg)^{n}-l \bigg)*x+ \\ 
& \bigg(\frac{1}{\sqrt{5}}\bigg(\frac{1}{2}+\frac{\sqrt{5}}{2}\bigg)^{n-1}-\frac{1}{\sqrt{5}}\bigg(\frac{1}{2}-\frac{\sqrt{5}}{2}\bigg)^{n-1}-s \bigg) 	\\	 
\end{flalign*}
Then $\phi^t=\tilde  x$ and so	$t=\log_\phi \tilde x$.



\section{Silver ratio}

\begin{figure}[H]
\centering
\begin{tikzpicture}

		\node[node] at (-1,1.5) {0-tile};
		\draw[axis] (0,1.5) -- (2*pi*2.2,1.5);
		\draw[shift={(0,1.5)}, thick] (0,3pt) -- (0,-3pt) node[below] {$0$};
		\path[shift={(0,1.5)}, thick](pi*2.2,3pt) -- (pi*2.2,0pt) node[above]{$b$};
		\path[shift={(0,1.5)}, thick](pi*2.2,3pt) -- (pi*2.2,-3pt) node[below]{$r$};
		\draw[shift={(0,1.5)}, thick] (2*pi*2.2,3pt) -- (2*pi*2.2,-3pt) node[below] {$1$};
		
		\node[node] at (-1,0) {1-tile};
		\draw[axis] (0,0) -- (2*pi*2.2,0);
		\draw[thick] (0,3pt) -- (0,-3pt) node[below] {$0$};
		\path[thick](pi*2.2,3pt) -- (pi*2.2,0pt) node[above]{$a$};
		\path[thick](pi*2.2,3pt) -- (pi*2.2,-3pt) node[below]{$1$};
		\draw[thick] (2*pi*2.2,3pt) -- (2*pi*2.2,-3pt) node[below] {$1$};
		
		\node[node] at (-1,-1.5) {2-tile};
		\draw[axis] (0,-1.5) -- (2*pi*2.2,-1.5);
	\draw[shift={(0,-1.5)}, thick] (0,3pt) -- (0,-3pt) node[below] {$0$};
		\path[shift={(0,-1.5)}] (0.4142*pi*2.2,3pt) -- (0.4142*pi*2.2,0pt) node[above] {$a$};
		\path[shift={(0,-1.5)}] (0.4142*pi*2.2,3pt) -- (0.4142*pi*2.2,-3pt) node[below] {$\frac{1}{r}$};
		\draw[shift={(0,-1.5)}, thick] (0.8284*pi*2.2,3pt) -- (0.8284*pi*2.2,-3pt) ;
		\path[shift={(0,-1.5)}] (1.24263*pi*2.2,3pt) -- (1.24263*pi*2.2,0pt) node[above] {$a$};
		\path[shift={(0,-1.5)}] (1.24263*pi*2.2,3pt) -- (1.24263*pi*2.2,-3pt) node[below] {$\frac{1}{r}$};
		\draw[shift={(0,-1.5)}, thick] (1.65684*pi*2.2,3pt) -- (1.65684*pi*2.2,-3pt) ;
		\path[shift={(0,-1.5)}] (1.82841*pi*2.2,3pt) -- (1.82841*pi*2.2,0pt) node[above] {$b$};
		\path[shift={(0,-1.5)}] (1.82841*pi*2.2,3pt) -- (1.82841*pi*2.2,-3pt) node[below] {$\frac{1}{r^2}$};
		\draw[shift={(0,-1.5)}, thick] (2*pi*2.2,3pt) -- (2*pi*2.2,-3pt) node[below] {$1$};

		\node[shift={(0,-3)}, node] at (-1,0) {3-tile};
		\draw[shift={(0,-3)}, axis] (0,0) -- (2*pi*2.2,0);     
		\draw[shift={(0,-3)}, thick] (0,3pt) -- (0,-3pt) node[below] {$0$};
		\path[shift={(0,-3)}] (0.17157*pi*2.2,3pt) -- (0.17157*pi*2.2,-3pt) node[below] {$\frac{1}{r^2}$};
		\path[shift={(0,-3)}] (0.17157*pi*2.2,3pt) -- (0.17157*pi*2.2,0pt) node[above] {$a$};
		\draw[shift={(0,-3)}, thick] (0.34314*pi*2.2,3pt) -- (0.34314*pi*2.2,-3pt);
		
		\path[shift={(0,-3)}] (0.51437*pi*2.2,3pt) -- (0.51437*pi*2.2,-3pt) node[below] {$\frac{1}{r^2}$};
		\path[shift={(0,-3)}] (0.51437*pi*2.2,3pt) -- (0.51437*pi*2.2,0pt) node[above] {$a$};
		\draw[shift={(0,-3)}, thick] (0.68628*pi*2.2,3pt) -- (0.68628*pi*2.2,-3pt);
		
		\path[shift={(0,-3)}] (0.75735*pi*2.2,3pt) -- (0.75735*pi*2.2,-3pt) node[below] {$\frac{1}{r^3}$};
		\path[shift={(0,-3)}] (0.75735*pi*2.2,3pt) -- (0.75735*pi*2.2,0pt) node[above] {$b$};
		\draw[shift={(0,-3)}, thick] (0.82842*pi*2.2,3pt) -- (0.82842*pi*2.2,-3pt);
		
		\path[shift={(0,-3)}] (1*pi*2.2,3pt) -- (1*pi*2.2,-3pt) node[below] {$\frac{1}{r^2}$};
		\path[shift={(0,-3)}] (1*pi*2.2,3pt) -- (1*pi*2.2,0pt) node[above] {$a$};
		\draw[shift={(0,-3)}, thick] (1.1716*pi*2.2,3pt) -- (1.1716*pi*2.2,-3pt);
			
		\path[shift={(0,-3)}] (1.3431*pi*2.2,3pt) -- (1.3431*pi*2.2,-3pt) node[below] {$\frac{1}{r^2}$};
		\path[shift={(0,-3)}] (1.3431*pi*2.2,3pt) -- (1.3431*pi*2.2,0pt) node[above] {$a$};
		\draw[shift={(0,-3)}, thick] (1.5147*pi*2.2,3pt) -- (1.5147*pi*2.2,-3pt);
		
		\path[shift={(0,-3)}] (1.5858*pi*2.2,3pt) -- (1.5858*pi*2.2,-3pt) node[below] {$\frac{1}{r^3}$};
		\path[shift={(0,-3)}] (1.5858*pi*2.2,3pt) -- (1.5858*pi*2.2,0pt) node[above] {$b$};
		\draw[shift={(0,-3)}, thick] (1.6569*pi*2.2,3pt) -- (1.6569*pi*2.2,-3pt);
		
		\path[shift={(0,-3)}] (1.8284*pi*2.2,3pt) -- (1.8284*pi*2.2,-3pt) node[below] {$\frac{1}{r^2}$};
		\path[shift={(0,-3)}] (1.8284*pi*2.2,3pt) -- (1.8284*pi*2.2,0pt) node[above] {$a$};
		
		\draw[shift={(0,-3)}, thick] (2*pi*2.2,3pt) -- (2*pi*2.2,-3pt) node[below] {$1$};

		\node[shift={(0,-4.5)}, node] at (-1,0) {4-tile};
		\draw[shift={(0,-4.5)}, axis] (0,0) -- (2*pi*2.2,0);
		\draw[shift={(0,-4.5)}, thick] (0,3pt) -- (0,-3pt) node[below] {$0$};
		
		\path[shift={(0,-4.5)}] (0.0711*pi*2.2,3pt) -- (0.0711*pi*2.2,-3pt) node[below] {$\frac{1}{r^3}$};
		\path[shift={(0,-4.5)}] (0.0711*pi*2.2,3pt) -- (0.0711*pi*2.2,0pt) node[above] {$a$};
		\draw[shift={(0,-4.5)}, thick] (0.1421*pi*2.2,3pt) -- (0.1421*pi*2.2,-3pt);
		
		\path[shift={(0,-4.5)}] (0.21321*pi*2.2,3pt) -- (0.21321*pi*2.2,-3pt) node[below] {$\frac{1}{r^3}$};
		\path[shift={(0,-4.5)}] (0.21321*pi*2.2,3pt) -- (0.21321*pi*2.2,0pt) node[above] {$a$};
		\draw[shift={(0,-4.5)}, thick] (0.28428*pi*2.2,3pt) -- (0.28428*pi*2.2,-3pt);
		
		\path[shift={(0,-4.5)}] (0.31372*pi*2.2,3pt) -- (0.31372*pi*2.2,-3pt) node[below] {$\frac{1}{r^4}$};
		\path[shift={(0,-4.5)}] (0.31372*pi*2.2,3pt) -- (0.31372*pi*2.2,0pt) node[above] {$b$};
		\draw[shift={(0,-4.5)}, thick] (0.34316*pi*2.2,3pt) -- (0.34316*pi*2.2,-3pt);
		
		\path[shift={(0,-4.5)}] (0.4142*pi*2.2,3pt) -- (0.4142*pi*2.2,-3pt) node[below] {$\frac{1}{r^3}$};
		\path[shift={(0,-4.5)}] (0.4142*pi*2.2,3pt) -- (0.4142*pi*2.2,0pt) node[above] {$a$};
		\draw[shift={(0,-4.5)}, thick] (0.4853*pi*2.2,3pt) -- (0.4853*pi*2.2,-3pt);
		
		\path[shift={(0,-4.5)}] (0.5563*pi*2.2,3pt) -- (0.5563*pi*2.2,-3pt) node[below] {$\frac{1}{r^3}$};
		\path[shift={(0,-4.5)}] (0.5563*pi*2.2,3pt) -- (0.5563*pi*2.2,0pt) node[above] {$a$};
		\draw[shift={(0,-4.5)}, thick] (0.6274*pi*2.2,3pt) -- (0.6274*pi*2.2,-3pt);
		
		\path[shift={(0,-4.5)}] (0.6569*pi*2.2,3pt) -- (0.6569*pi*2.2,-3pt) node[below] {$\frac{1}{r^4}$};
		\path[shift={(0,-4.5)}] (0.6569*pi*2.2,3pt) -- (0.6569*pi*2.2,0pt) node[above] {$b$};
		\draw[shift={(0,-4.5)}, thick] (0.6863*pi*2.2,3pt) -- (0.6863*pi*2.2,-3pt);
		
		\path[shift={(0,-4.5)}] (0.7574*pi*2.2,3pt) -- (0.7574*pi*2.2,-3pt) node[below] {$\frac{1}{r^3}$};
		\path[shift={(0,-4.5)}] (0.7574*pi*2.2,3pt) -- (0.7574*pi*2.2,0pt) node[above] {$a$};
		\draw[shift={(0,-4.5)}, thick] (0.8284*pi*2.2,3pt) -- (0.8284*pi*2.2,-3pt);
		
		\path[shift={(0,-4.5)}] (0.8995*pi*2.2,3pt) -- (0.8995*pi*2.2,-3pt) node[below] {$\frac{1}{r^3}$};
		\path[shift={(0,-4.5)}] (0.8995*pi*2.2,3pt) -- (0.8995*pi*2.2,0pt) node[above] {$a$};
		\draw[shift={(0,-4.5)}, thick] (0.9706*pi*2.2,3pt) -- (0.9706*pi*2.2,-3pt);
		
		\path[shift={(0,-4.5)}] (1.04167*pi*2.2,3pt) -- (1.04167*pi*2.2,-3pt) node[below] {$\frac{1}{r^3}$};
		\path[shift={(0,-4.5)}] (1.04167*pi*2.2,3pt) -- (1.04167*pi*2.2,0pt) node[above] {$a$};
		\draw[shift={(0,-4.5)}, thick] (1.11274*pi*2.2,3pt) -- (1.11274*pi*2.2,-3pt);
		
		\path[shift={(0,-4.5)}] (1.14218*pi*2.2,3pt) -- (1.14218*pi*2.2,-3pt) node[below] {$\frac{1}{r^4}$};
		\path[shift={(0,-4.5)}] (1.14218*pi*2.2,3pt) -- (1.14218*pi*2.2,0pt) node[above] {$b$};
		\draw[shift={(0,-4.5)}, thick] (1.1716*pi*2.2,3pt) -- (1.1716*pi*2.2,-3pt);
		
		\path[shift={(0,-4.5)}] (1.2426*pi*2.2,3pt) -- (1.2426*pi*2.2,-3pt) node[below] {$\frac{1}{r^3}$};
		\path[shift={(0,-4.5)}] (1.2426*pi*2.2,3pt) -- (1.2426*pi*2.2,0pt) node[above] {$a$};
		\draw[shift={(0,-4.5)}, thick] (1.3137*pi*2.2,3pt) -- (1.3137*pi*2.2,-3pt);
		
		\path[shift={(0,-4.5)}] (1.38483*pi*2.2,3pt) -- (1.38483*pi*2.2,-3pt) node[below] {$\frac{1}{r^3}$};
		\path[shift={(0,-4.5)}] (1.38483*pi*2.2,3pt) -- (1.38483*pi*2.2,0pt) node[above] {$a$};
		\draw[shift={(0,-4.5)}, thick] (1.4559*pi*2.2,3pt) -- (1.4559*pi*2.2,-3pt);
		
		\path[shift={(0,-4.5)}] (1.48534*pi*2.2,3pt) -- (1.48534*pi*2.2,-3pt) node[below] {$\frac{1}{r^4}$};
		\path[shift={(0,-4.5)}] (1.48534*pi*2.2,3pt) -- (1.48534*pi*2.2,0pt) node[above] {$b$};
		\draw[shift={(0,-4.5)}, thick] (1.5147*pi*2.2,3pt) -- (1.5147*pi*2.2,-3pt);
		
		\path[shift={(0,-4.5)}] (1.5858*pi*2.2,3pt) -- (1.5858*pi*2.2,-3pt) node[below] {$\frac{1}{r^3}$};
		\path[shift={(0,-4.5)}] (1.5858*pi*2.2,3pt) -- (1.5858*pi*2.2,0pt) node[above] {$a$};
		\draw[shift={(0,-4.5)}, thick] (1.6569*pi*2.2,3pt) -- (1.6569*pi*2.2,-3pt);
		
		\path[shift={(0,-4.5)}] (1.7279*pi*2.2,3pt) -- (1.7279*pi*2.2,-3pt) node[below] {$\frac{1}{r^3}$};
		\path[shift={(0,-4.5)}] (1.7279*pi*2.2,3pt) -- (1.7279*pi*2.2,0pt) node[above] {$a$};
		\draw[shift={(0,-4.5)}, thick] (1.799*pi*2.2,3pt) -- (1.799*pi*2.2,-3pt);
		
		\path[shift={(0,-4.5)}] (1.87013*pi*2.2,3pt) -- (1.87013*pi*2.2,-3pt) node[below] {$\frac{1}{r^3}$};
		\path[shift={(0,-4.5)}] (1.87013*pi*2.2,3pt) -- (1.87013*pi*2.2,0pt) node[above] {$a$};
		\draw[shift={(0,-4.5)}, thick] (1.9412*pi*2.2,3pt) -- (1.9412*pi*2.2,-3pt);
		
		\path[shift={(0,-4.5)}] (1.97064*pi*2.2,3pt) -- (1.97064*pi*2.2,-3pt) node[below] {$\frac{1}{r^4}$};
		\path[shift={(0,-4.5)}] (1.97064*pi*2.2,3pt) -- (1.97064*pi*2.2,0pt) node[above] {$b$};
		
		\draw[shift={(0,-4.5)}, thick] (2*pi*2.2,3pt) -- (2*pi*2.2,-3pt) node[below] {$1$};	
\end{tikzpicture}
\caption{Silver ratio tilings}
\label{silvertil}
\end{figure}
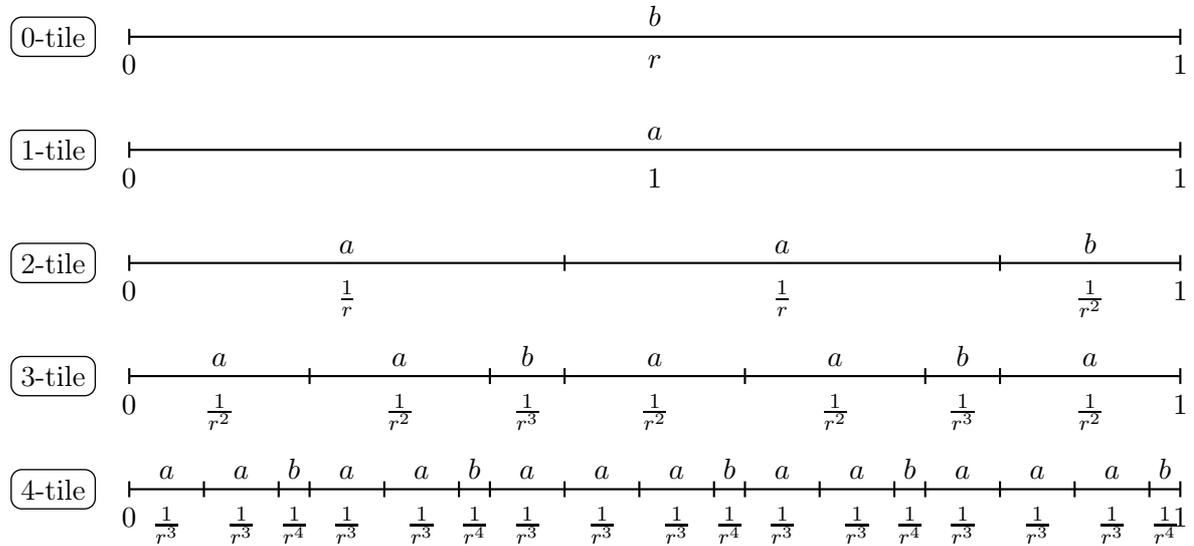
\subsection{Silver ratio basics}
In this section we will demonstrate another example using a different tiling of $[0,1]$, the tiling corresponding to the silver ratio. Recall the definition of the silver ratio. It is the positive root of $r^2-2r-1=0 $. The silver ratio is denoted by $\delta=1+\sqrt 2$.

Recall also that the silver ratio corresponds to a recurrence relation, similarly to how the golden ratio corresponds to the Fibonacci recurrence relation. 

The silver ratio recurrence relation is $a_{n}=2a_{n-1}+a_{n-2}$. The silver ratio  sequence is 0, 1, 2, 5, 12, 29, 70, 169, ... where we used initial conditions $a_{0}=0$ and $a_{1}=1$. The closed form formula for the $a_n$ is $a_{n}=\frac{1}{2\sqrt{2}}\big(1+ \sqrt{2} \big)^{n}-\frac{1}{2\sqrt{2}}\big(1- \sqrt{2}	\big)^n$. \\ 

The silver ratio also corresponds to a two-character sequence as follows. We again begin with a word with a single letter $a$, and apply a subtitution $b\to a$, $a\to aab$, Applying the substitution repeatedly, we obtain

$$a\to aab\to aabaaba\to aabaabaaabaabaaab\to\ldots$$
. We can then obtain a sequence of tilings of $[0,1]$ similarly as before. 
\subsection{Fractal Dimension for Silver ratio fractals}

\begin{figure}[H]
\centering
\begin{tikzpicture}

		\node[node] at (-1.5,-1.5) {tiling};
		\draw[axis] (0,-1.5) -- (2*pi,-1.5);
	\draw[shift={(0,-1.5)}, thick] (0,3pt) -- (0,-3pt) node[below] {$0$};
		\path[shift={(0,-1.5)}] (0.4142*pi,3pt) -- (0.4142*pi,0pt) node[above] {$a$};
		\path[shift={(0,-1.5)}] (0.4142*pi,3pt) -- (0.4142*pi,-3pt) node[below] {$\frac{1}{\delta}$};
		\draw[shift={(0,-1.5)}, thick] (0.8284*pi,3pt) -- (0.8284*pi,-3pt) ;
		\path[shift={(0,-1.5)}] (1.24263*pi,3pt) -- (1.24263*pi,0pt) node[above] {$a$};
		\path[shift={(0,-1.5)}] (1.24263*pi,3pt) -- (1.24263*pi,-3pt) node[below] {$\frac{1}{\delta}$};
		\draw[shift={(0,-1.5)}, thick] (1.65684*pi,3pt) -- (1.65684*pi,-3pt) ;
		\path[shift={(0,-1.5)}] (1.82841*pi,3pt) -- (1.82841*pi,0pt) node[above] {$b$};
		\path[shift={(0,-1.5)}] (1.82841*pi,3pt) -- (1.82841*pi,-3pt) node[below] {$\frac{1}{\delta^2}$};
		\draw[shift={(0,-1.5)}, thick] (2*pi,3pt) -- (2*pi,-3pt) node[below] {$1$};
		
	\draw[arrow] (pi,-2.7) --(pi,-3.1);

		\node[shift={(-1,-3.8)}, node] at (-1,0) {1st removal};
		
		\draw[shift={(0,-3.8)}, axis] (0,0) -- (0.8284*pi,0);
		\draw[shift={(0,-3.8)}, axis] (1.6569*pi,0) -- (2*pi,0);
		
		\draw[shift={(0,-3.8)}, thick] (0,3pt) -- (0,-3pt) node[below] {$0$};
		\path[shift={(0,-3.8)}] (0.4142*pi,3pt) -- (0.4142*pi,0pt) node[above] {$a$};
		\path[shift={(0,-3.8)}] (0.4142*pi,3pt) -- (0.4142*pi,-3pt) node[below] {$\frac{1}{\delta}$};

		\draw[shift={(0,-3.8)}, thick] (0.82842*pi,3pt) -- (0.82842*pi,-3pt);

		\draw[shift={(0,-3.8)}, thick] (1.6569*pi,3pt) -- (1.6569*pi,-3pt);
		
		\path[shift={(0,-3.8)}] (1.8284*pi,3pt) -- (1.8284*pi,-3pt) node[below] {$\frac{1}{\delta^2}$};
		\path[shift={(0,-3.8)}] (1.8284*pi,3pt) -- (1.8284*pi,0pt) node[above] {$b$};
		
		\draw[shift={(0,-3.8)}, thick] (2*pi,3pt) -- (2*pi,-3pt) node[below] {$1$};

	\draw[arrow] (pi,-5) --(pi,-5.4);

		\node[shift={(-1,-6.1)}, node] at (-1,0) {2nd removal};
		
		\draw[shift={(0,-6.1)}, axis] (0,0) -- (0.34314*pi,0);
		\draw[shift={(0,-6.1)}, axis] (0.68628*pi,0) -- (0.82842*pi,0);
		\draw[shift={(0,-6.1)}, axis] (1.6569*pi,0) -- (1.798974*pi,0);
		\draw[shift={(0,-6.1)}, axis] (1.941108*pi,0) -- (2*pi,0);
		
		\draw[shift={(0,-6.1)}, thick] (0,3pt) -- (0,-3pt) node[below] {$0$};
		\draw[shift={(0,-6.1)}, thick] (0.34314*pi,3pt) -- (0.34314*pi,-3pt);
	
		\draw[shift={(0,-6.1)}, thick] (0.68628*pi,3pt) -- (0.68628*pi,-3pt);
		
		\draw[shift={(0,-6.1)}, thick] (0.82842*pi,3pt) -- (0.82842*pi,-3pt);

		\draw[shift={(0,-6.1)}, thick] (1.6569*pi,3pt) -- (1.6569*pi,-3pt);
		\draw[shift={(0,-6.1)}, thick] (1.798974*pi,3pt) -- (1.798974*pi,-3pt);
		
		\draw[shift={(0,-6.1)}, thick] (1.941108*pi,3pt) -- (1.941108*pi,-3pt);
		
		\draw[shift={(0,-6.1)}, thick] (2*pi,3pt) -- (2*pi,-3pt) node[below] {$1$};
		
	\end{tikzpicture}
\caption{The construction of a (2,1,0) fractal.}
\label{s2t}
\end{figure}
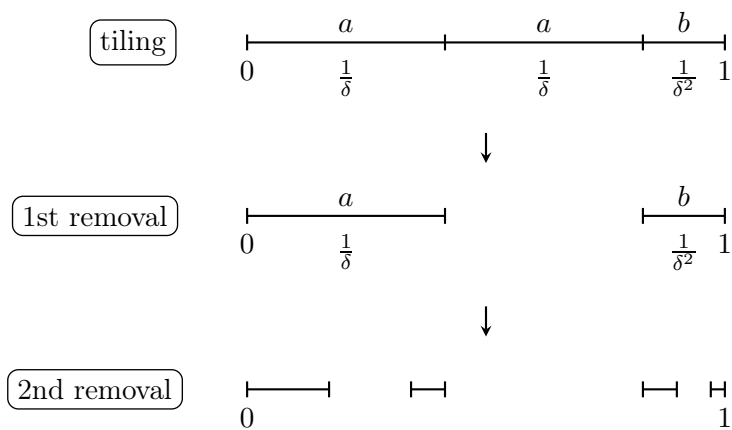
Let us consider the $(2,1,0)$ fractal corresponding to this silver ratio tiling. Recall that this means we look at the word generated by applying the $b\to a, a\to aab$ substitution twice to $b$ (thus obtaining the word $aab$) and then generate the tiling of $[0,1]$ corresponding to this word, which consists now of two ``long tiles" of length $1/\delta$, and one ``short tile" of length $1/\delta^2$. We then remove $1$ long tile, and zero short tiles. 

We then perform the same process iteratively on each of the remaining intervals, in a manner similar to the construction of the Cantor set.

The similarity dimension, $d$ of the resulting (2,1,0) fractal is 
\begin{flalign*}
{(scaling\ factor, r)}^{(dimension, d)}   &= length\ of\ 2-tile'  \\
(\delta^d)					    &=\frac{1}{\delta^2}+(\frac{1}{\delta} *\frac{1}{\delta^d}) \\
	\delta^{2d}			 &=\delta^d+1 \\
	let\ \delta^d=x &\\
x^{2}				    &=x+1 \\
0&=x^2-x-1 \\
x&\approx 1.6180 \\
so\ \delta^d \approx 1.6180 \\
d					        &= \log_{\delta}1.6180 \\
						    & \approx 0.54596
\end{flalign*}

The Hausdorff dimension is calculated similarly as before:
\begin{flalign*}
(\frac{1}{\delta})^t+(\frac{1}{\delta^2})^t &=1 \\
\delta^t+1 &= \delta^2t \\
	0    &=\delta^{2t}-\delta^t-1 \\
	let\ \delta^t=x &\\
x^{2}				    &=x+1 \\
0&=x^2-x-1 \\
x&\approx 1.6180 \\
so\ \delta^d \approx 1.6180 \\
t					        &= \log_{\delta}1.6180 \\
						    & \approx 0.54596
\end{flalign*}

The general similarity dimension for the $(n,l,s)$ fractal corresponding to the silver ratio is then calculated as follows:

\begin{flalign*}
{(scaling\ factor, r)}^{(dimension, d)}   =&length\ after\ scale,\ l_a\ \\
(\delta^d)						=&(total\ length\ of\ long\ tile, t_l\ )+(total\ length\ of\ short\ tile, t_s\ )*\delta^{-d} \\
							=&(a_{n}-l)+(a_{n-1}-s)*\delta^{-d} \\
(\delta^d)^{n}					=&(a_{n}-l)*\delta^d+(a_{n-1}-s) \\
Let\ \tilde{x}\ be\ the\ solution\ to\  x^n  =&\bigg(\frac{1}{2\sqrt{2}}\bigg(1+ \sqrt{2} \bigg)^{n}-\frac{1}{2\sqrt{2}}\bigg(1- \sqrt{2}	\bigg)^n \bigg)-l \bigg)*x+ \\ 
&\bigg(\frac{1}{2\sqrt{2}}\bigg(1+ \sqrt{2} \bigg)^{n-1}-\frac{1}{2\sqrt{2}}\bigg(1- \sqrt{2}	\bigg)^{n-1} \bigg)-s \bigg) 	\\	 
Then\ \delta^d=& x\\					
d=&\log_\delta x
\end{flalign*}

For Hausdorff dimension in a general $(n,l,s)$ fractal:
\begin{flalign*}
H^t(x)=&(number\ of\ intervals)(the\ length\ of\ the\ intervals)^t \\
1=&(total\ length\ of\ long\ tiles\ intervals, t_l)+ \\
&(total\ length\ of\ the\ short\ tile\ intervals, t_s) \\
1=&(a_{n}-l)(1/\delta^3)^t+(a_{n-1}-s)(1/\delta^4)^t \\
1=&(1/\delta^t)^4[(a_{n}-l)\delta^t+(a_{n-1}-s)] \\
\delta^{4t}=&(a_{n}-l)\delta^t+(a_{n-1}-s) \\
(\delta^t)^n=&(a_{n}-l)\delta^t+(a_{n-1}-s) \\
Let\ \tilde{x}\ be\ the\ solution\ to\  x^n  =&\bigg(\frac{1}{2\sqrt{2}}\bigg(1+ \sqrt{2} \bigg)^{n}-\frac{1}{2\sqrt{2}}\bigg(1- \sqrt{2}	\bigg)^n \bigg)-l \bigg)*x+ \\ 
&\bigg(\frac{1}{2\sqrt{2}}\bigg(1+ \sqrt{2} \bigg)^{n-1}-\frac{1}{2\sqrt{2}}\bigg(1- \sqrt{2}	\bigg)^{n-1} \bigg)-s \bigg) 	\\	 
Then\ \delta^t=& x\\					
t=&\log_\delta x
\end{flalign*}


\section{Fractals corresponding to general metallic means sequences}
Let us now consider fractals corresponding to general metallic means sequences. Our substitution is now
\begin{flalign*}
\rho_{p,q} : \ \ \ \ \ a			&\rightarrow a^{p}b^{q} \\
b			&\rightarrow a 
\end{flalign*}
\\
where $p$ and $q$ are positive integers.
The $p,q$ metallic mean is the positive solution to the polynomial $r^2=pr+q$, that is $r=(-p+\sqrt{p^2+4q})/2$. 

This $r$ will again serve as the scaling factor for our tiling of $[0,1]$.  

We can then generate $(n,l,s)$ fractals as before: we consider the tiling that corresponds to applying the substitution to $b$ $n$ times, remove $l$ long tiles and $s$ short tiles, and then repeat the process iteratively. 

The process to calculate the similarity and Hausdorff dimensions works the same way as before:
\begin{flalign*}
r^d  =&length\ after\ substitution,l_a\  \\
(\gamma^d)				=&(total\ length\ of\ long\ tile,t_l\ )+(total\ length\ of\ short\ tile,t_s\ )*\gamma^{-d} \\
							=&(a_{n}-l)+((a_{n-1}-s)*\gamma^{-d}) \\
(\gamma^{n})^{d}				=&(a_{n}-l)*\gamma^{d}+(a_{n-1}-s) 
\end{flalign*}
Let $\tilde{x}$ be the positive solution to  
\begin{flalign*}
 x^{n}=& \bigg(\frac{1}{\sqrt{p^{2}+4q}}\bigg(\frac{1}{2}\bigg(\sqrt{p^{2}+4q}+p\bigg)\bigg)^n  -\frac{1}{\sqrt{p^{2}+4q}}\bigg(\frac{1}{2}\bigg(p-\sqrt{p^{2}+4q}\bigg)\bigg)^n-l\bigg)x+  \\
& \bigg(\frac{1}{\sqrt{p^{2}+4q}}\bigg(\frac{1}{2}\bigg(\sqrt{p^{2}+4q}+p\bigg)\bigg)^{n-1}-\frac{1}
{\sqrt{p^{2}+4q}}\bigg(\frac{1}{2}\bigg(p-\sqrt{p^{2}+4q}\bigg)\bigg)^{n-1}-s\bigg) 
\end{flalign*}

Then $\gamma^{d}				=\tilde x$, and so 
$d=\log_{\gamma}{\tilde x}$.

The tilings were covered by the number of intervals, $i$ with length of the intervals, $j$ powered by $t$, so the general Hausdorff dimension of the substituted p,q-tilings is shown as

For Hausdorff dimension in a general $(n,l,s)$ fractal:
\begin{flalign*}
H^t(x)=&(number\ of\ intervals)(the\ length\ of\ the\ intervals)^t \\
1=&(total\ length\ of\ long\ tiles\ intervals, t_l)+ \\
&(total\ length\ of\ the\ short\ tile\ intervals, t_s) \\
1=&(a_{n}-l)(1/\gamma^{n-1})^t+(a_{n-1}-s)(1/\gamma^n)^t \\
1=&(1/\gamma^t)^n[(a_{n}-l)\gamma^t+(a_{n-1}-s)] \\
\gamma^{nt}=&(a_{n}-l)\gamma^t+(a_{n-1}-s) \\
(\gamma^t)^n=&(a_{n}-l)\gamma^t+(a_{n-1}-s) \\
Let\ \tilde{x}\ & be\ the\ positive\ solution\ to\  \\
x^n  =& \bigg(\frac{1}{\sqrt{p^{2}+4q}}\bigg(\frac{1}{2}\bigg(\sqrt{p^{2}+4q}+p\bigg)\bigg)^n  -\frac{1}{\sqrt{p^{2}+4q}}\bigg(\frac{1}{2}\bigg(p-\sqrt{p^{2}+4q}\bigg)\bigg)^n-l\bigg)x+  \\
& \bigg(\frac{1}{\sqrt{p^{2}+4q}}\bigg(\frac{1}{2}\bigg(\sqrt{p^{2}+4q}+p\bigg)\bigg)^{n-1}-\frac{1}
{\sqrt{p^{2}+4q}}\bigg(\frac{1}{2}\bigg(p-\sqrt{p^{2}+4q}\bigg)\bigg)^{n-1}-s\bigg) 	\\	
Then\ \gamma^t=& x\\					
t=&\log_\gamma x
\end{flalign*}

\begin{paragraph}{Acknowledgements}
D.~C.~O.\ and Y.~S.~J.~L\ were supported in part by a grant from the Fundamental Research Grant Scheme from the Malaysian Ministry of Education (grant number FRGS/1/2022/TK07/XMU/01/1). D.~C.~O.\ was also supported by a grant from the National Natural Science Foundation of China  (grant number 12201524), and a Xiamen University Malaysia Research Fund (grant number XMUMRF/2023C11/IMAT/0024).
\end{paragraph}

\bibliographystyle{alpha}
\bibliography{abcd1234}

\end{document}